\documentclass[11pt]{article}
\usepackage[a4paper,margin=1in]{geometry}
\usepackage[authoryear,round]{natbib}
\usepackage[T1]{fontenc}
\usepackage{lmodern,microtype}
\usepackage{amsmath,amssymb,amsthm,mathtools,mathrsfs}
\usepackage{booktabs,array,enumitem,needspace,graphicx}
\usepackage{xcolor}
\usepackage[colorlinks=true,linkcolor=blue!45!black,citecolor=blue!45!black,urlcolor=blue!45!black]{hyperref}
\hypersetup{pdftitle={When are dividend barriers optimal? A distributional characterization in renewal risk models},pdfauthor={Linyu Miao}}
\numberwithin{equation}{section}
\newtheorem{theorem}{Theorem}[section]
\newtheorem{lemma}[theorem]{Lemma}
\newtheorem{proposition}[theorem]{Proposition}
\newtheorem{corollary}[theorem]{Corollary}
\theoremstyle{definition}
\newtheorem{definition}[theorem]{Definition}
\newtheorem{assumption}[theorem]{Assumption}

\newcommand{\E}{\mathbb E}
\newcommand{\Pp}{\mathbb P}

\newcommand{\ind}{\mathbf 1}
\newcommand{\J}{\mathcal J}

\newcommand{\Sop}{\mathcal S}
\newcommand{\Var}{\operatorname{Var}}
\newcommand{\Exp}{\operatorname{Exp}}

\newcommand{\suppref}[1]{Appendix~\ref{#1}}
\newcommand{\mainref}[1]{\ref{#1}}
\newcommand{\mainequation}[1]{\eqref{#1}}

\begin{document}
\title{When are dividend barriers optimal?\\
A distributional characterization in renewal risk models}
\author{Linyu Miao\thanks{Corresponding author. Email: \texttt{mly202210@126.com}.}\\
\small School of Mathematical Sciences, Nankai University\\
\small Tianjin, 300071, China}
\date{October 30, 2026}
\maketitle

\begin{abstract}
We study the finite-horizon optimal dividend problem in a renewal
risk model with unrestricted dividend payments, examining how
claim-size distributions determine barrier optimality. Under general
assumptions on a fixed interclaim-time distribution, we prove that
exponential claims
admit a globally optimal barrier for every positive safety loading,
positive discount rate and finite horizon. Conversely, within the
class of claim-size distributions with continuous nondecreasing hazard,
every nonexponential law admits one fixed premium--discount pair for
which no barrier is globally optimal at any sufficiently long horizon.
These results establish, within this distribution class, an equivalence
between exponentially distributed claims and the existence of a globally
optimal barrier for every parameter combination in the stated range.
Numerical examples compare exponential and Gamma claims with the same
mean and illustrate how the claim-size distribution changes the
structure of dividend decisions.
\end{abstract}

\section{Introduction}
The optimal dividend problem is a central problem in actuarial risk
theory: how much of an insurer's surplus should be distributed to
shareholders, and how much should be retained to meet future claims?
The criterion introduced by \citet{DeFinetti1957} evaluates a strategy
through the expected present value of dividends paid before ruin.
An immediate payment reduces the funds available to support future
payments; see
\citet{Avanzi2009,Schmidli2008}. Barrier strategies are among the most
widely studied forms of dividend policy. Their rule is transparent:
pay any surplus above a prescribed level and retain the remainder.
This simplicity makes it important to understand when the rule is
optimal, and how that conclusion depends on the claim-size distribution.

Classical results establish barrier optimality for several important
risk models. \citet{AsmussenTaksar1997} obtain an optimal barrier in
a diffusion model with unrestricted dividends.
\citet{AvramPalmowskiPistorius2007} study the dividend problem for
spectrally negative L\'evy processes and express barrier values through
scale functions. \citet{Loeffen2008} gives a sufficient condition in
terms of complete monotonicity of the jump density, while
\citet{KyprianouRiveroSong2010,LoeffenRenaud2010} relate barrier
optimality to convexity properties of the jump distribution and the
associated scale function. These results provide substantial classes
of models in which a barrier is optimal. They do not make the barrier
form universal: the examples of \citet{AzcueMuler2005} and the
computational study of \citet{AlbrecherGarcia2023} also show the
relevance of dividend strategies with several bands.

Renewal arrivals introduce a further consideration. The distribution
of the next claim time depends on how long the insurer has already
waited since the last claim. \citet{AlbrecherClaramuntMarmol2005}
study dividend payments under barriers with generalized Erlang
interclaim times, and \citet{MishuraSchmidli2012} investigate barrier
strategies in the same renewal setting.
\citet{AlbrecherHartinger2006} show by a counterexample that a constant
barrier need not be optimal in a Sparre Andersen model. This leaves
open an important distinction: failure of a constant level does not
exclude a barrier that adjusts to the time elapsed since the last claim.
A finite planning horizon also allows the barrier to vary with time; \citet{DeAngelisEkstrom2017} establish such a time-dependent
barrier in a diffusion dividend problem.

For general renewal models, \citet{TianBaiGuo2020} establish properties
of the singular-dividend value function and characterize it through
a constrained viscosity formulation of the Hamilton--Jacobi--Bellman
equation. \citet{LiuLiuLiu2020} study bounded dividend rates and obtain
verification and a space--time band structure under their conditions.
With exponentially distributed claims and observable phases of a
phase-type waiting time, \citet{TianLiu2022} study optimal phase-dependent
barriers. These contributions address the value equation or optimal
policies in specified settings. The question considered here is which
claim-size laws guarantee a barrier across premium rates, discount
rates and horizons, when the controller observes claim times and sizes
but no latent phases.

We answer this question for each fixed renewal law with continuous
hazard and uniformly bounded expected residual waiting time. Our first
main result, Theorem~\ref{thm:classification}, proves that, within the
class of claim-size distributions with continuous nondecreasing hazard,
a globally optimal barrier exists for every positive safety loading,
positive discount rate and finite horizon \emph{if and only if the
claim-size distribution is exponential}. The barrier may depend on
time and claim age, meaning the time since the last claim.
For exponential claims we construct this barrier, prove its continuity,
and establish concavity and continuous differentiability of the value
in surplus. The same distributional characterization holds for
age-dependent barriers on the infinite horizon.

The second result, Theorem~\ref{thm:conflict}, identifies a stronger
failure for every nonexponential law in this class. There is one fixed
premium--discount pair for which \emph{no globally optimal barrier
exists at any sufficiently long finite horizon}, or on the infinite
horizon. At the same time and claim age, all of a small
initial surplus should be paid immediately, whereas full payment is
strictly suboptimal at a larger surplus. These decisions cannot be
implemented by a single upper barrier. They also imply that the value
is not concave in surplus. Thus extending the horizon cannot restore
barrier optimality at these parameters, even when the barrier is allowed
to depend on both time and the observed claim history through age.

The main difficulty is that optimal payments depend not only on
current surplus but also on the time since the last claim and the
terminal date. The associated HJB equation involves three variables---
time, surplus and claim age---and an integral term linking the current
state to post-claim states corresponding to different claim sizes.
This dependence on several variables and the nonlocal term make the
optimal payment region and its boundary difficult to determine directly.
We separate payments before the next claim from decisions afterwards.
For exponential claims, the memoryless property permits a comparison of the benefits
of retaining different amounts of capital. We show that the resulting
payment decisions fit together into one barrier and verify its
optimality. For nonexponential claims, we choose premium income close
to the expected claim outgo and compare the consequences of paying
now or retaining capital. At small reserves, delaying payment creates
a loss that an optimal strategy cannot sustain; at a larger reserve,
a short delay yields a strictly higher value. These comparisons prove
the failure of a barrier without solving for the full optimal strategy
or assuming that the value is differentiable.

Our numerical examples compare exponential and Gamma claims while
keeping the mean claim size, the claim-arrival mechanism and the
economic parameters fixed. For exponential claims, the computed
payment decisions have a single-barrier structure. For Gamma claims,
full payment at small surplus coexists with capital retention at a
larger surplus, at the same time and claim age. These decisions cannot
be implemented by one upper barrier. This comparison is consistent
with the theoretical analysis and shows how differences in the shape
of the claim-size distribution can change dividend decisions even
when the mean claim size is unchanged.

Section~\ref{sec:model} introduces the risk model and states the two
main results. Section~\ref{sec:dpp} establishes dynamic programming and
basic value bounds. Section~\ref{sec:positive} constructs the optimal
barrier for exponential claims, and Section~\ref{sec:negative} proves
its failure for nonexponential claims at suitable fixed parameters.
Section~\ref{sec:examples} discusses claim distributions and numerical
examples. Auxiliary estimates are collected in the appendices;
further examples and numerical details appear in
Appendices~\ref{supp:numerics}--\ref{app:certificate}.

\section{Risk model and main results}\label{sec:model}
\subsection{The renewal risk model}
Let $(\Omega,\mathcal F,\mathbb P)$ be a complete probability space
supporting two independent sequences $(A_i)_{i\ge1}$ and
$(Y_i)_{i\ge1}$ of independent identically distributed positive random
variables. The variables $A_i$ and $Y_i$ represent interclaim times and
claim sizes. Denote their distribution functions by $F_A$ and $G$,
their densities by $f_A$ and $g$, and set
\[
\overline F_A=1-F_A,\qquad S=1-G,\qquad
m=\E A_1\in(0,\infty),\qquad \mu=\E Y_1\in(0,\infty).
\]
We assume $\overline F_A(a)>0$ for finite $a$ and that the renewal
hazard $\eta=f_A/\overline F_A$ is continuous on $[0,\infty)$.

For an initial time $t\ge0$ and claim age $a\ge0$, let $A_a$ have survival
\[
S_a(s):=\mathbb P(A_a>s)
=\frac{\overline F_A(a+s)}{\overline F_A(a)},\qquad s\ge0.
\]
Thus $A_a$ is the additional waiting time for the next claim, given
that the current waiting period has already lasted $a$ units;
$S_a(s)$ is the probability of waiting at least another $s$ units.
Under the delayed renewal law $\mathbb P_a$, the first claim occurs
at $T_1=t+A_a$, subsequent interclaim times have law $F_A$, and all
claim sizes have law $G$, independently. For $s\ge t$, write
\[
T_n=t+A_a+\sum_{j=2}^n A_j,\qquad
N_s^{t,a}=\sum_{n\ge1}\ind_{\{T_n\le s\}},\qquad
C_s=\sum_{j=1}^{N_s^{t,a}}Y_j .
\]
Here $N_s^{t,a}$ counts claims between the initial time $t$ and time
$s$, and $C_s$ is their total amount, with $C_{t-}=0$.
The filtration $\mathbb F^{t,a}=(\mathcal F_s^{t,a})_{s\ge t}$ is the
completed, right-continuous natural filtration of this marked claim
process. The claim-age process is
\[
a_s=
\begin{cases}
a+s-t,&N_s^{t,a}=0,\\
s-T_{N_s^{t,a}},&N_s^{t,a}\ge1.
\end{cases}
\]
Thus the controller observes claim times and sizes, and claim age
resets to zero at each claim. The conditional law of future arrivals
depends on this observed age, even when a phase-type representation
of the waiting-time law is available.

For premium rate $c>0$ and initial surplus $x\ge0$, the uncontrolled
reserve at time $s\ge t$ is $x+c(s-t)-C_s$. If $L_s$ is the cumulative
dividend paid from time $t$ to time $s$, with $L_{t-}=0$, the controlled
reserve before killing is
\begin{equation*}
X_s^L=x+c(s-t)-C_s-L_s,\qquad s\ge t,\qquad X_{t-}^L=x.
\end{equation*}
The main results use $cm>\mu$, the usual net-profit condition.

\begin{assumption}\label{ass:arrival}
The interclaim-time law satisfies
\begin{equation*}
M:=\sup_{a\ge0}m_A(a)<\infty,\qquad
m_A(a):=\E A_a=\int_0^\infty S_a(s)\,ds .
\end{equation*}
\end{assumption}
The quantity $m_A(a)$ is the expected additional wait at age $a$;
the assumption bounds this expectation uniformly over the elapsed wait.
Together with the preceding regularity of $\eta$, this is the standing
arrival assumption for the classification. The uniform bound is used
in the converse; the finite-horizon exponential construction needs
only continuous hazard and positive survival.

\subsection{Dividend strategies and the value function}
Fix a maturity $T>0$, an initial time $t\in[0,T]$ and a discount
rate $q>0$. A dividend strategy is an $\mathbb F^{t,a}$-adapted,
nondecreasing, right-continuous process $(L_s)_{s\ge t}$, locally of
finite variation, with $L_{t-}=0$. An initial lump sum
$L_t\in[0,x]$ is allowed. Dividend payments cannot cause ruin.
More precisely, at a claim the claim is settled before any dividend:
with $\Delta C_s=C_s-C_{s-}$,
\[
\widehat X_s=X_{s-}^L-\Delta C_s,\qquad
0\le\Delta L_s\le\widehat X_s
\quad\hbox{when }\widehat X_s\ge0 .
\]
If $\widehat X_s<0$, ruin occurs and $\Delta L_s=0$.
The continuous part of $L$ must also preserve nonnegative surplus
before ruin. Define the ruin time by
\[
\tau^L=\inf\{s\ge t:\widehat X_s<0\},\qquad
\inf\varnothing=\infty .
\]
The strategy is stopped at $T\wedge\tau^L$. The set of such admissible
strategies is denoted by $\mathcal A^T(t,x,a)$. At maturity, any
surviving reserve may be liquidated. Zero surplus is not absorbing:
premium income can be paid at rate $c$ while waiting for the next claim.

Write $\E_{t,x,a}$ for expectation from the initial state $(t,x,a)$.
The performance criterion and value function are
\begin{align*}
J^T(t,x,a;L)
&=\E_{t,x,a}\!\left[
 \int_{[t,T\wedge\tau^L]}e^{-q(s-t)}\,dL_s\right],\\
U^T(t,x,a)&=\sup_{L\in\mathcal A^T(t,x,a)}J^T(t,x,a;L).
\end{align*}
The integral includes the initial payment at $t$ and any terminal
payment at $T$. There is no payment at ruin, since
$\Delta L_{\tau^L}=0$ on $\{\tau^L<\infty\}$. On the infinite horizon, time homogeneity allows
us to take the initial time to be zero. Denote the corresponding
admissible class by $\mathcal A^\infty(x,a)$ and write
$\E_{x,a}=\E_{0,x,a}$. Then
\begin{align*}
J^\infty(x,a;L)&=\E_{x,a}\!\left[
 \int_{[0,\tau^L]}e^{-qs}\,dL_s\right],\\
V(x,a)&=\sup_{L\in\mathcal A^\infty(x,a)}J^\infty(x,a;L),
\end{align*}
where the integral is understood as an increasing limit if
$\tau^L=\infty$. In particular,
\begin{equation}\label{eq:terminal}
U^T(T,x,a)=x .
\end{equation}
For fixed $T$ we write $U=U^T$. The delayed formulation is defined on
$[0,T]\times[0,\infty)^2$; states of a process started at time zero
also satisfy $a\le t$. The counterexample states below satisfy this
physical restriction.

\subsection{The HJB equation and barrier strategies}
For a smooth test function $f$, the killed, uncontrolled generator is
\[
\begin{aligned}
\mathcal L f(t,x,a)
&=c f_x(t,x,a)+f_a(t,x,a)\\
&\quad+\eta(a)\left[\int_0^x f(t,x-y,0)g(y)\,dy-f(t,x,a)\right].
\end{aligned}
\]
The claim integral includes surviving claims; the continuation after
ruin is zero. Immediate payment of $d\in[0,x]$ gives
\[
U(t,x,a)\ge d+U(t,x-d,a).
\]
Thus the dynamic programming equation is the nonlocal HJB variational
inequality
\begin{equation}\label{eq:HJB}
\max\left\{U_t+\mathcal L U-qU,\ 1-U_x\right\}=0,
\qquad 0\le t<T,\quad x>0,\quad a\ge0,
\end{equation}
with terminal condition \eqref{eq:terminal}.
At $x=0$ the dynamics are subject to the solvency constraint;
there is no absorbing boundary condition $U(t,0,a)=0$.
Where one-sided derivatives exist, the boundary equation is
\[
\max_{0\le d\le c}
\{U_t+U_a+(c-d)U_x-(q+\eta(a))U+d\}(t,0,a)=0.
\]
On the infinite horizon \eqref{eq:HJB} has no $U_t$ term and $U$ is
replaced by $V$. These equations express the local optimality
conditions. The arguments below use the integral dynamic programming
principle, so no global classical regularity of $U$ is assumed.

\begin{definition}\label{def:barrier}
A time--age barrier for maturity $T$ is a measurable function
$b_T:[0,T]\times[0,\infty)\to[0,\infty)$ whose initial projection
and minimal upper regulation define an admissible strategy.
Writing $b=b_T$, the initial dividend from $(t,x,a)$ is
$(x-b(t,a))^+$. Before the first claim its regulator is
\begin{equation}\label{eq:reflection}
L_s^b=\sup_{t\le r\le s}[x+c(r-t)-b(r,a+r-t)]^+,
\qquad t\le s\le T .
\end{equation}
After each surviving claim the rule restarts at the new time, surplus and age zero. It is globally optimal if this one
feedback attains $U^T$ from every initial state. On the physical
domain the requirement is restricted to $a\le t$.
An infinite-horizon age barrier $b_\infty(a)$ is defined in the
same way.
\end{definition}
The definition includes nonsmooth barriers whenever the regulator is
admissible; continuity is proved for the exponential optimizer.

\subsection{Main results}
\begin{definition}\label{def:ifr}
Let $\mathcal I$ be the class of claim laws with finite positive mean,
$S(x)>0$ for finite $x$, and continuous, nondecreasing hazard
$h_Y(x)=g(x)/S(x)$ on $[0,\infty)$. These laws have an increasing failure rate (IFR), meaning a
nondecreasing hazard. No global upper bound on the hazard is imposed.
\end{definition}
The first result asks for one optimal barrier for each choice of
parameters, with that barrier serving every initial state. Throughout
the following statements, $b_T$ and $b_\infty$ range over the admissible
barriers of Definition~\ref{def:barrier}.

\begin{theorem}[Characterization by barrier optimality]\label{thm:classification}
Fix an interclaim-time law satisfying Assumption~\ref{ass:arrival}
and the regularity conditions above, and let $G\in\mathcal I$.
The following are equivalent:
\begin{enumerate}[label=(\roman*),leftmargin=*]
\item $G=\Exp(\beta)$ for some $\beta>0$;
\item For every $c>\mu/m$, $q>0$ and $T>0$, there is a barrier
$b_T$ such that, for all $(t,x,a)\in[0,T]\times[0,\infty)^2$,
\[
J^T(t,x,a;L^{b_T})=U^T(t,x,a).
\]
\item For every $c>\mu/m$ and $q>0$, there is an age-dependent
barrier $b_\infty$ such that, for all $(x,a)\in[0,\infty)^2$,
\[
J^\infty(x,a;L^{b_\infty})=V(x,a).
\]
\end{enumerate}
\end{theorem}
Failure of this guarantee alone need only occur at one choice of
parameters and maturity. The next result gives a stronger conclusion:
one fixed premium--discount pair excludes barriers at every
sufficiently long maturity.

\begin{theorem}[Persistent failure of barrier optimality]\label{thm:conflict}
Under the assumptions of Theorem~\ref{thm:classification}, the
following are equivalent:
\begin{enumerate}[label=(\roman*),leftmargin=*]
\item $G$ is nonexponential;
\item There exist $c>\mu/m$, $q>0$ and $T_0>0$ such that,
for every $T\ge T_0$ and every barrier $b_T$, some initial state
$(t,x,a)\in[0,T]\times[0,\infty)^2$ satisfies
\[
J^T(t,x,a;L^{b_T})<U^T(t,x,a).
\]
\end{enumerate}
When these conditions hold, the parameters can be chosen together
with $a_*\le t_*<T_0$ and $0<x_*<x^*$ such that, for every $T\ge T_0$,
\[
U^T(t_*,x,a_*)=x+U^T(t_*,0,a_*),\qquad 0\le x\le x_*.
\]
For $0<x\le x_*$ and $L\in\mathcal A^T(t_*,x,a_*)$,
\begin{equation}\label{eq:main-small}
J^T(t_*,x,a_*;L)=U^T(t_*,x,a_*)\quad\Longrightarrow\quad L_{t_*}=x.
\end{equation}
At the larger surplus,
\begin{equation}\label{eq:main-large}
U^T(t_*,x^*,a_*)>x^*+U^T(t_*,0,a_*).
\end{equation}
The same assertions hold for $V$, with initial payment $L_0$ and
the same $c,q,a_*,x_*,x^*$.
In particular, for every age-dependent barrier $b_\infty$, some
$(x,a)\in[0,\infty)^2$ satisfies
$J^\infty(x,a;L^{b_\infty})<V(x,a)$.
Both indicated surplus sections of the value are nonconcave.
\end{theorem}
The small-surplus conclusion concerns the initial payment; the insurer
can still receive premiums and pay further dividends afterwards.
At the larger surplus the theorem excludes full payment, while allowing
partial payment followed by retention of a positive reserve.

\section{Dynamic programming and value bounds}\label{sec:dpp}
Before the next claim, no new random observation arrives. We can
therefore describe the initial part of a strategy by its payments
along the no-claim path, leaving the post-claim decisions in the value
function. This representation will be used both to construct a barrier
and to compare competing payment decisions. In this representation,
$s\ge0$ denotes time elapsed after $t$: on the no-claim path, the
payment plan satisfies $\ell_s=L_{t+s}$, so $\ell_0=L_t$.

For a function $f(v,x,a)$, put
\[
\J_g f(v,z,0)=\int_0^z f(v,z-y,0)g(y)\,dy,\qquad
w_a(s)=e^{-qs}S_a(s).
\]
The first quantity is the expected continuation after a claim from
reserve $z$, assigning zero to ruin. The weight $w_a(s)$ discounts a
payment at elapsed time $s$ and accounts for the probability that no
claim has yet occurred. Use $\J_\beta$ when $g(y)=\beta e^{-\beta y}$.

\begin{definition}[Claim-free payment plan]\label{def:plan}
For $0\le t<T$ and $x\ge0$, a feasible claim-free payment plan is a
nondecreasing right-continuous function $\ell:[0,T-t]\to[0,\infty)$
satisfying
\[
\ell_{0-}=0,\qquad \ell_s\le x+cs\quad(0\le s\le T-t),\qquad
\ell_{T-t}=x+c(T-t).
\]
Its retained surplus is $z_s=x+cs-\ell_s$.
\end{definition}
The final equality requires all remaining reserve to be paid at
maturity on the no-claim path:
\[
\Delta\ell_{T-t}=x+c(T-t)-\ell_{(T-t)-},\qquad z_{T-t}=0.
\]
This restriction does not change the value: terminal dividends are
allowed, and unpaid terminal reserve earns no further reward. At
$t=T$ the plan consists only of the payment $x$.

The following proposition collects the dynamic programming identity
and the bounds needed below. It also supplies optimal deterministic
plans at fixed finite-horizon states; the optimal barrier itself will
be constructed and verified in Section~\ref{sec:positive}.
\begin{proposition}\label{prop:dpp}
The value $U$ is continuous and satisfies
\begin{equation}\label{eq:firstclaim}
U(t,x,a)=x+c\int_0^{T-t} w_a(s)\,ds
+\sup_z\int_0^{T-t} w_a(s)H_U(t,a;s,z_s)\,ds,
\end{equation}
where the supremum is over the plans of Definition~\ref{def:plan} and
\begin{equation}\label{eq:H}
H_f(t,a;s,z)=\eta(a+s)\J_gf(t+s,z,0)-(q+\eta(a+s))z.
\end{equation}
For every fixed finite-horizon state the supremum is attained.
Moreover,
\begin{align}
x+c\int_0^{T-t}w_a(s)\,ds&\le U(t,x,a)
\le x+\frac cq(1-e^{-q(T-t)}),\label{eq:budget}\\
0\le U(t,x,a)-U(t+\delta,x,a)&\le c\delta,\quad 0\le\delta\le T-t,\nonumber\\
0\le V(x,a)-U(t,x,a)&\le \frac cq e^{-q(T-t)}.
\label{eq:tails}
\end{align}
The infinite-horizon counterpart of \eqref{eq:firstclaim} holds with
$V$, an infinite integral, and nondecreasing right-continuous plans
satisfying $\ell_{0-}=0$ and $\ell_s\le x+cs$ for all $s\ge0$.
\end{proposition}
\begin{proof}
At a fixed initial state, conditioning on the first claim gives the
objective
\[
\Phi_f(\ell)=\int_{[0,T-t]}w_a(s)\,d\ell_s+
\int_0^{T-t} w_a(s)\eta(a+s)\J_gf(t+s,z_s,0)\,ds.
\]

Since $d\ell_s=c\,ds-dz_s$, $z_{0-}=x$, $z_{T-t}=0$ and
$w_a'=-(q+\eta(a+\cdot))w_a$, integration by parts gives
\[
\begin{aligned}
\int_{[0,T-t]}w_a(s)\,d\ell_s
&=c\int_0^{T-t}w_a(s)\,ds
  -[w_a(s)z_s]_{0-}^{T-t}
  +\int_0^{T-t}w_a'(s)z_s\,ds\\
&=x+c\int_0^{T-t}w_a(s)\,ds
  -\int_0^{T-t}w_a(s)(q+\eta(a+s))z_s\,ds .
\end{aligned}
\]
This gives the objective in \eqref{eq:firstclaim}; the dynamic
programming justification follows below. The residual waiting time
has a density. Conditional on the observed history at a claim epoch,
a no-claim plan is deterministic; its countably many payment discontinuities therefore
have zero probability of coinciding with the next claim. Pre-payment
and post-payment versions at these times give the same continuation
integral.

For completeness, let $U_n$ be the value terminated at the $n$th future
claim, with no liquidation at that terminating claim, and let $U_0=0$.
For $n\ge1$ liquidation at maturity remains allowed. Write $\Sop$ for the first-claim operator.
Lemma~\ref{lem:plans} in Appendix~\ref{app:plans} supplies a fixed
metrizable plan space, a continuous compact-valued feasible
correspondence, and a continuous integrated objective.
The lemma gives continuity and a Borel maximizing plan at every
finite step. Concatenating these plans through the finitely many
remaining claims proves $U_n=\Sop U_{n-1}$. At $t=T$, $U_n(T,x,a)=x$ for $n\ge1$;
the separate convention $U_0=0$ describes termination before any payment.

The pathwise total dividends up to maturity are at most $x+c(T-t)$.
Terminating a full policy at the $n$th claim thus loses at most
\begin{equation}\label{eq:claimtail}
0\le U-U_n
\le(x+c(T-t))\Pp(N_T^{t,a}\ge n)
\le(x+c(T-t))\Pp(N_{T-t}^{0,0}\ge n-1).
\end{equation}

For $\lambda>0$ and $n\ge2$,
\[
\mathbb P(N_{T-t}^{0,0}\ge n-1)
=\mathbb P\!\left(\sum_{j=1}^{n-1}A_j\le T-t\right)
\le e^{\lambda(T-t)}(\E e^{-\lambda A_1})^{n-1}.
\]
Since $\E e^{-\lambda A_1}<1$, the convergence is uniform for bounded
$T-t$. Hence $U_n\to U$ locally uniformly and $U$ is continuous. On compact state sets the first-claim operator is
continuous in the locally uniform norm, since its source is evaluated
only at surplus at most $x+c(T-t)$. Thus $U=\Sop U$.
Lemma~\ref{lem:plans}, applied to the continuous source $U$, gives
attainment for each fixed state.

We will also use continuation within any prescribed error $\delta>0$.
After a first-claim plan from $(t,x,a)$, surviving continuation states
have surplus at most $x+c(T-t)$ and remaining time at most $T-t$.
Choose $n$ so that \eqref{eq:claimtail} is at most $\delta$ throughout
this set. The Borel optimal policies of the truncated problem then
provide a continuation whose conditional value is at least $U-\delta$.
Thus any fixed plan with first-claim objective $\Phi_U(\ell)$ can be
extended to an admissible policy $L^\delta$ with
\begin{equation}\label{eq:near-continuation}
J^T(t,x,a;L^\delta)\ge\Phi_U(\ell)-\delta.
\end{equation}
The same construction at the first claim or a fixed elapsed time
$h\le T-t$ proves the corresponding stopped dynamic programming
identity.

Immediate distribution followed by payment of premium income proves
the lower bound in \eqref{eq:budget}; the discounted capital budget
proves the upper bound. To shorten a policy by $\delta$, liquidate its
remaining surplus at the earlier maturity. Subsequent continuation
is bounded by that surplus plus at most $c\delta$ of additional premium
income. 
For $0\le\delta\le T-t$, time homogeneity and this truncation give
\[
U(t+\delta,x,a)\le U(t,x,a)\le U(t+\delta,x,a)+c\delta.
\]
For an infinite policy, liquidation at elapsed time $T-t$ removes
at most the additional discounted premium budget:
\[
J^\infty(x,a;L)
\le U(t,x,a)+e^{-q(T-t)}\frac cq .
\]
Taking suprema proves \eqref{eq:tails}.

The second bound in \eqref{eq:tails} makes $V$ continuous as a
locally uniform limit. To pass the first-claim identity to the infinite
horizon, use the raw dividend and continuation objective above. For
$0\le f(z,0)\le z+c/q$, its part after elapsed time $H$, uniformly
over feasible plans and bounded initial surplus, is bounded by
\begin{equation}\label{eq:plan-tail}
C(1+H)e^{-qH}.
\end{equation}
Indeed, the dividend tail is at most $e^{-qH}(x+cH+c/q)$, and
\[
\int_H^\infty w_a(s)\eta(a+s)\J_g f(z_s,0)\,ds
\le\sup_{s\ge H}e^{-qs}(x+cs+c/q),
\]
since $\int_H^\infty S_a(s)\eta(a+s)\,ds\le1$.
On $[0,H]$, \eqref{eq:tails} gives uniform convergence of the
continuation sources as maturity tends to infinity. Any infinite
plan can be liquidated at a finite maturity, and any finite plan can
be extended by paying all subsequent premiums; their initial parts
agree. First let maturity increase with $H$ fixed, then let
$H\to\infty$. The uniform tail bound \eqref{eq:plan-tail} permits taking suprema and
proves the infinite-horizon identity. Integration by parts is valid
there because $0\le w_a(s)z_s\le e^{-qs}(x+cs)\to0$.
\end{proof}

The bound \eqref{eq:tails} quantifies the cost of imposing a maturity:
after the remaining reserve is liquidated, only future premium income
can be lost. Its discounted value is at most $(c/q)e^{-q(T-t)}$.

\Needspace{14\baselineskip}
\section{Exponential claims: construction of an optimal barrier}
\label{sec:positive}
Throughout this section $G=\Exp(\beta)$. The finite-horizon arguments
require only positive survival and a continuous interclaim hazard.
The positive result holds for every $c,q>0$, so this section does
not require a positive safety loading.

\subsection{The value of additional surplus}
An extra unit of capital can be paid immediately, but retaining it
may also preserve future dividend opportunities after a claim. The
next estimate bounds this additional benefit. Set
\[
r_a(t)=\E[e^{-qA_a};A_a\le T-t]<1.
\]
This is the expected present value of one unit received at the next
claim if it occurs before maturity, and zero otherwise. It combines
the chance of a claim with discounting; $q>0$ and $A_a>0$ almost surely
make the weight strictly less than one.

\begin{lemma}\label{lem:expbudget}
For $n\ge1$ and $0\le x<y$,
\begin{equation}\label{eq:expbudget}
1\le\frac{U_n(t,y,a)-U_n(t,x,a)}{y-x}
\le1+r_a(t)[\beta U_{n-1}(t,0,0)-1]^+.
\end{equation}
The full value satisfies the same inequality with $U$ on both
occurrences. On the infinite horizon, replace $U$ by $V$ and
$r_a(t)$ by $\E e^{-qA_a}<1$.
\end{lemma}
\begin{proof}
Immediate payment of $y-x$ gives the lower bound. For the upper bound,
first put $y=x+d$ and copy an upper-state policy at the lower state,
clipping a payment if the available lower surplus is insufficient.
Let $R_s$ be the cumulative excess paid by the upper process by time
$s$, stopped before lower ruin or termination, with $R_{t-}=0$; while both copies are active, their
surplus gap is $\delta_s=d-R_s\in[0,d]$. The copied upper history is
computable from the common observations, so clipping defines an
admissible lower policy.

At a claim time, condition on the pre-claim history. If the lower
pre-claim surplus is $z$, then, conditional on $Y>z$, the overshoot
$Y-z$ has law $\Exp(\beta)$. With gap $\delta_-$, the probability
of upper survival is $1-e^{-\beta\delta_-}\le\beta\delta_-$,

and its surviving surplus belongs to $[0,d]$. Conditional on lower
ruin at time $s$, the surviving upper continuation is bounded by
\[
\int_0^{\delta_-}\beta e^{-\beta u}
 U_{n-1}(s,\delta_--u,0)\,du
\le \beta\delta_- U_{n-1}(t,d,0).
\]
Here surplus monotonicity and the time inequality in
\eqref{eq:tails}, valid also for $U_{n-1}$, give the last step.
This continuation includes post-claim dividends, which are therefore
excluded from $R$ at the stopping claim.

Let $\sigma$ be lower ruin, the terminating claim or maturity,
whichever comes first. Suppress terminal payments in $R$ and
write $\delta_{\sigma-}$ for the remaining pre-settlement gap.
Stieltjes integration by parts gives
\[
\int_{[t,\sigma)}e^{-q(s-t)}\,dR_s
+e^{-q(\sigma-t)}\delta_{\sigma-}
=d-q\int_t^\sigma e^{-q(s-t)}\delta_s\,ds\le d.
\]
At maturity the payoff difference left after the pre-terminal
payments is this gap; at the terminating claim it is zero.
At a ruin claim, the conditional estimate above bounds its expected
replacement by $\beta\delta_{\sigma-}U_{n-1}(t,d,0)$.
Consequently, writing $J^+$ and $J^-$ for the two expected payoffs,
\[
J^+-J^-\le d+
d[\beta U_{n-1}(t,d,0)-1]^+
\E[e^{-q(\sigma-t)};\sigma\le T,\ \text{lower ruin at a claim}].
\]
The last expectation is at most $r_a(t)$ because $\sigma-t\ge A_a$
on that event. Taking near-optimal upper policies gives this bound
for the value difference. Subdivide $[x,y]$ into $k$ increments
$d=(y-x)/k$ and let $k\to\infty$; continuity at zero proves
\eqref{eq:expbudget}. The same argument without claim truncation
uses $U$ in the continuation. On the infinite horizon the surviving
gap is bounded by $d e^{-qs}$ and vanishes as $s\to\infty$.
\end{proof}

The capital estimate has a structural consequence: the value is
concave in surplus. The accompanying convolution property is what
will allow the optimal payment decisions at different capital levels
to be assembled into one barrier.
\begin{lemma}\label{lem:expconcavity}
For every $n\ge1$, $U_n(t,\cdot,a)$ is concave. If
$f(x)=U_n(v,x,0)$ and $B=f(0)$, then:
\begin{enumerate}[label=(\roman*),leftmargin=*]
\item $\beta B\le1$ implies $f(x)=x+B$;
\item $\beta B>1$ implies $f'_+(0)<\beta B$ and $\J_\beta f$
is strictly concave.
\end{enumerate}
These assertions also hold for $U$ and $V$.
\end{lemma}
\begin{proof}
$U_1=x+c\int_0^{T-t} w_a$ is affine. Suppose the assertions of
concavity have been proved at step $n-1$, and put
$f=U_{n-1}(v,\cdot,0)$, $B=f(0)$. If $\beta B\le1$,
\eqref{eq:expbudget} and $U_{n-2}\le U_{n-1}$ show that $f=x+B$.
Thus
\[
(\J_\beta f)'(x)=1+(\beta B-1)e^{-\beta x}\le1,
\]
and $H_{U_{n-1}}(t,a;s,\cdot)$ is strictly decreasing. For $n=1$ the source is zero
and the same conclusion is immediate.

If $\beta B>1$, \eqref{eq:expbudget} and $r_0(v)<1$ give
$f'_+(0)<\beta B$. Integrating the derivative measure of $f'_+$ gives,
for almost every $x>0$,
\begin{equation}\label{eq:convcurv}
(\J_\beta f)''(x)
=\beta e^{-\beta x}\{f'_+(0)-\beta B\}
-\beta\int_{(0,x]}e^{-\beta(x-y)}(-df'_+)(dy)<0.
\end{equation}
This identity follows by integrating the exponential convolution of
$f'_+$ once more by parts, and needs no classical second derivative of $f$.

The function $B=U_{n-1}(v,0,0)$ is nonincreasing in $v$.
For fixed $(t,a)$, the times at which
$\beta U_{n-1}(t+s,0,0)\le1$ consequently form a terminal interval.
On that interval $H_{U_{n-1}}(t,a;s,z)$ is strictly decreasing from zero.

Let $s_*$ be the beginning of that interval, with $s_*=T-t$ if the
interval reduces to maturity. Replacing a feasible retained path by
$z_s=0$ for $s\ge s_*$ weakly increases its objective. For
$s<s_*$ the source is concave. Choose maximizing retained paths
$z^1,z^2$ from $x_1,x_2$, zero after $s_*$, and let
\[
x_\lambda=\lambda x_1+(1-\lambda)x_2,\qquad
z^\lambda=\lambda z^1+(1-\lambda)z^2,\qquad 0\le\lambda\le1.
\]
The corresponding payments
$\ell^\lambda=\lambda\ell^1+(1-\lambda)\ell^2$ are feasible from
$x_\lambda$. Therefore, writing $H=H_{U_{n-1}}$,
\[
\begin{aligned}
U_n(t,x_\lambda,a)
&\ge x_\lambda+c\int_0^{T-t}w_a(s)\,ds
 +\int_0^{s_*}w_a(s)H(t,a;s,z_s^\lambda)\,ds\\
&\ge \lambda U_n(t,x_1,a)+(1-\lambda)U_n(t,x_2,a).
\end{aligned}
\]
This proves the induction step.
The two assertions at the new step follow from
Lemma~\ref{lem:expbudget} and \eqref{eq:convcurv}.
The locally uniform limit \eqref{eq:claimtail} preserves concavity.
Then \eqref{eq:tails} and the infinite-horizon version of
Lemma~\ref{lem:expbudget} prove the assertions for $V$.
\end{proof}

\subsection{Construction and verification of the barrier}
We now turn concavity into a control. The next representation separates
the contributions of initial capital and future premiums, and determines
when each should be paid. For the full value let
\begin{equation}\label{eq:marginal}
\begin{aligned}
F(v,z)&=\J_\beta U(v,z,0),\\
\psi(v,z,\alpha)&=\eta(\alpha)\{F_z(v,z)-1\}-q,\\
I(t,x,a;h)&=\int_0^h w_a(s)\psi(t+s,x+cs,a+s)\,ds,\\
\mathcal M(t,x,a)&=\max_{0\le h\le T-t}I(t,x,a;h).
\end{aligned}
\end{equation}
The identity $F_z=\beta(U-F)$ ensures joint continuity without
presupposing differentiability of $U$.

\begin{lemma}[Representation of the value]\label{lem:layers}
The value has a continuous spatial derivative, with a right derivative
at zero, and
\begin{align}
U_x(t,x,a)&=1+\mathcal M(t,x,a),\label{eq:derivative}\\
U(t,x,a)&=x+
c\int_0^{T-t} w_a(s)\{1+\mathcal M(t+s,0,a+s)\}\,ds
+\int_0^x\mathcal M(t,y,a)\,dy.\label{eq:layerformula}
\end{align}
\end{lemma}
\begin{proof}
Write a feasible path as $z_s=cs+y_s$, where $y_s=x-\ell_s$ is
nonincreasing. Label surplus layers by $y\in[-c(T-t),x]$.
The birth time of layer $y$ is $s_y=\max(0,-y/c)$: nonnegative labels
are initial surplus, whereas negative labels represent future premiums.
Its payment time $\vartheta_y$ belongs to $[s_y,T-t]$.
Feasibility requires $\vartheta_y$ to be nonincreasing in $y$:
higher layers must be paid no later than lower ones.
Fubini's theorem and $H_U(t,a;s,0)=0$ give
\begin{equation}\label{eq:layerdecomp}
\int_0^{T-t} w_a(s)H_U(t,a;s,z_s)\,ds
=\int_{-c(T-t)}^x
\int_{s_y}^{\vartheta_y}
w_a(s)\psi(t+s,y+cs,a+s)\,ds\,dy.
\end{equation}

For each $y$ define
\[
\vartheta_y^*=
\min\operatorname*{arg\,max}_{s_y\le h\le T-t}
\int_{s_y}^h w_a(s)\psi(t+s,y+cs,a+s)\,ds .
\]
We show that these choices satisfy the required ordering.
Indeed, on the terminal interval identified in
Lemma~\ref{lem:expconcavity}, $\psi\le-q$, so earlier-born layers stop
no later than its beginning and later-born layers stop immediately.
Before that interval the integrand is nonincreasing in $y$.
If $y_1<y_2$ but $\vartheta_{y_1}^*<\vartheta_{y_2}^*$, both layers can
stop at these two times. Optimality implies
\[
\begin{aligned}
0
&\le\int_{\vartheta_{y_1}^*}^{\vartheta_{y_2}^*}
 w_a(s)\psi(t+s,y_2+cs,a+s)\,ds\\
&\le\int_{\vartheta_{y_1}^*}^{\vartheta_{y_2}^*}
 w_a(s)\psi(t+s,y_1+cs,a+s)\,ds
\le0 .
\end{aligned}
\]
The integrand ordering forces equality, making
$\vartheta_{y_1}^*$ an earlier maximizer for layer $y_2$, a contradiction.
The earliest maximizers are jointly Borel measurable in the state and
layer label: the birth and terminal times are continuous, the
maximizer correspondence is measurable and compact-valued, and
minimizing the time coordinate over it gives its earliest member.
Define explicitly
\begin{equation}\label{eq:layer-control}
\ell_s^*=\int_{-c(T-t)}^{x}
\ind_{\{\vartheta_y^*\le s\}}\,dy,\qquad0\le s\le T-t.
\end{equation}
This is nondecreasing and right-continuous, by dominated convergence
from above in $s$. Only labels $y\ge-cs$ can have
$\vartheta_y^*\le s$, so $\ell_s^*\le x+cs$; also
$\ell_{T-t}^*=x+c(T-t)$. An atom at zero is allowed.
Ordering makes the paid labels an upper interval.
The unspent labels born by $s$ consequently form
$[-cs,x-\ell_s^*]$, up to endpoints of zero Lebesgue measure.
Their length is the retained surplus, so \eqref{eq:layerdecomp}
holds for this plan. It attains the integral of individual maxima
and is jointly measurable in the initial state.

For initial layers the maximum is $\mathcal M(t,y,a)$.
For premium layers substitute $y=-cr$ and use
$w_a(r+s)=w_a(r)w_{a+r}(s)$: their maximum is
$w_a(r)\mathcal M(t+r,0,a+r)$. Inserting these maxima into
\eqref{eq:firstclaim} gives \eqref{eq:layerformula}.
Continuity of $\mathcal M$ and differentiation of the last integral give
\eqref{eq:derivative}.
\end{proof}

\Needspace{12\baselineskip}
The preceding representation produces a feasible maximizing plan.
We next show that its payments follow a single continuous boundary,
so that the same rule can be used after every surviving claim.
\begin{theorem}[Exponential barrier]\label{thm:positive}
There is a finite continuous barrier
\begin{equation}\label{eq:barrier}
b(t,a)=\inf\{x\ge0:\mathcal M(t,x,a)=0\}
\end{equation}
that is globally optimal. Equivalently,
\begin{equation*}
x\ge b(t,a)\quad\Longleftrightarrow\quad
I(t,x,a;h)\le0\quad\text{for every }h\in[0,T-t].
\end{equation*}
Also $U_x>1$ below the barrier and $U_x=1$ above it.
For $T-t\le1/(\beta c)$ the barrier is zero at every age.
\end{theorem}
\begin{proof}
Layer maxima can ignore the terminal interval where $\psi<0$.
On its complement $\psi$ is nonincreasing in surplus, hence $\mathcal M$ is
nonincreasing. It is strictly decreasing wherever positive:
a maximizing positive integral must include positive claim intensity
in the interval of strictly concave $F$, where $F_z$ is strictly
decreasing.

To show finiteness, write $U(v,x,0)=x+R_v(x)$.
By concavity, $R_v$ is nonnegative, nondecreasing, concave and
bounded by $C=c(T-t)$ for $t\le v\le T$. Consequently
$R'_{v,+}(z)\le C/z$. Splitting the exponential convolution at
$x/2$ gives the uniform estimate
\begin{equation}\label{eq:exp-tail}
F_z(v,x)-1
\le\beta C(e^{-\beta x}+e^{-\beta x/2})+\frac{2C}{x}.
\end{equation}
Local boundedness of $\eta$ on $[a,a+T-t]$ now makes every integrand
in \eqref{eq:marginal} negative for sufficiently large $x$. Thus the
zero set of $\mathcal M$ is a nonempty upper interval.

Joint continuity of $\mathcal M$ gives lower semicontinuity of $b$.
Fix $x>b(t,a)$ and choose $b(t,a)<y<x$.
Every nonzero prefix $I(t,x,a;h)$ is strictly negative:
before the terminal interval use strict decrease of $F_z$ whenever
the intensity is positive; if no such times occur, the integral of
$-q$ is negative. After the terminal interval starts, prefixes only
decrease. Its initial integrand is also strictly negative, since
the initial integrand at $y$ is nonpositive and strict decrease applies,
or else $\psi\le-q$. Choose $h_1,\delta>0$ such that
$I(t,x,a;h)\le-\delta h$ for $0<h\le h_1$.
For $h\in[h_1,T-t]$ the continuous function $I$ has a strictly
negative maximum. Continuity of the integrand preserves the first
bound, and compactness preserves the second, under small changes
of $(t,a)$, including the change of the terminal endpoint.
Hence nearby parameters still have $\mathcal M(x)=0$, proving upper
semicontinuity of $b$. The zero terminal strip below handles $t=T$.

For a layer present at the initial time, let $\vartheta$ be its earliest
maximizing payment time. The identity
\begin{equation}\label{eq:prefix-flow}
I(t,x,a;r+h)=I(t,x,a;r)+
w_a(r)I(t+r,x+cr,a+r;h)
\end{equation}
shows that its remaining marginal maximum is positive before
$\vartheta$ and zero at $\vartheta$. Before the earliest global
maximizer the current prefix is strictly less than that maximum;
at the maximizer no subsequent prefix can improve it.
Thus $\vartheta$ is its first contact with \eqref{eq:barrier}.
The same identity holds for every premium layer.
A layer $y$ has been paid by $s$ precisely when
$y+cr\ge b(t+r,a+r)$ for some $r\le s$ after its birth.
Since $b\ge0$, this inequality itself excludes times before birth.
The paid labels therefore form the upper interval with lower endpoint
$\inf_{r\le s}[b(t+r,a+r)-cr]$.
Its length in \eqref{eq:layer-control} is exactly
$\sup_{0\le r\le s}[x+cr-b(t+r,a+r)]^+=L^b_{t+s}$,
the regulator in \eqref{eq:reflection} evaluated at time $t+s$.
The barrier is continuous and nonnegative, so this regulator is
right-continuous, nondecreasing and cannot overdraw surplus.
Restarting after surviving claims gives one admissible feedback strategy for
all states. Iterating the attained first-claim equality along this
specific strategy leaves a nonnegative remainder bounded by
\[
0\le U(t,x,a)-J^T(t,x,a;L^b)
\le (x+c(T-t))\mathbb P(N_T^{t,a}\ge n)\longrightarrow0.
\]
Thus $L^b$ is optimal over the full admissible class.

Finally $\beta U(v,0,0)\le\beta c(T-v)\le1$ on the stated
terminal strip. There $H_U$ is strictly decreasing throughout the
remaining interval, so immediate full distribution and subsequent premium
payments are optimal, giving $b=0$.
\end{proof}

\subsection{The infinite horizon}
The finite-horizon values converge to $V$, which preserves their
concavity. We construct the infinite-horizon barrier directly from $V$.
\begin{proposition}\label{prop:infinite-positive}
Exponential claims admit a globally optimal finite age barrier for
every $c,q>0$. With $F_\infty=\J_\beta V$ it is characterized by
\begin{align}
I_\infty(x,a;h)
&=\int_0^h w_a(s)
[\eta(a+s)\{F_\infty'(x+cs)-1\}-q]\,ds,\label{eq:Iinfty}\\
b_\infty(a)&=\inf\{x\ge0:\sup_{0\le h\le\infty}I_\infty(x,a;h)=0\}.
\nonumber
\end{align}
\end{proposition}
\begin{proof}
By \eqref{eq:tails}, finite-horizon concavity passes to $V$ and the
convolution derivatives converge, using
$(\J_\beta f)'=\beta(f-\J_\beta f)$.
If $\beta V(0,0)\le1$, Lemma~\ref{lem:expbudget} makes $V$ affine
in surplus; its claim source has derivative at most one, and the
zero barrier is optimal. Otherwise \eqref{eq:convcurv} makes
$F_\infty$ strictly concave with $F_\infty'(x)\downarrow1$.
Here $V(x,0)-x$ is bounded by $c/q$, and \eqref{eq:exp-tail} applies
with this bound.

The integrals in \eqref{eq:Iinfty} are continuous on the compactified
time interval $[0,\infty]$. In fact their tails are controlled by
\begin{equation}\label{eq:inftytail}
\int_H^\infty w_a(s)\eta(a+s)\,ds\le e^{-qH},
\qquad q\int_H^\infty w_a(s)\,ds\le e^{-qH},
\end{equation}
and $F_\infty'$ is bounded. Set
$C_F=\sup_{z\ge0}[F_\infty'(z)-1]^+<\infty$.
Writing $\mathcal M_\infty(x,a)=\max_h I_\infty(x,a;h)$,
the corresponding layer maxima satisfy
$0\le \mathcal M_\infty(x,a)\le C_F$, uniformly in age.
In addition to the single-layer bound \eqref{eq:inftytail}, future
premium layers born after $H$ have total excess value at most
\begin{equation*}
c\int_H^\infty w_a(r)\mathcal M_\infty(0,a+r)\,dr
\le(cC_F/q)e^{-qH}.
\end{equation*}
Thus Fubini and passage from finite birth intervals to all future
layers are justified. Their baseline premium value has the same
tail estimate with $C_F$ replaced by one. Formula
\eqref{eq:layer-control} extends to labels $y\in(-\infty,x]$:
on every bounded time interval only $[-cs,x]$ can be paid, so the
payment process is locally finite and right-continuous.
Earliest maximizing times in $[s_y,\infty]$ are measurable after
compactifying time, have the same order, and the prefix identity
\eqref{eq:prefix-flow} identifies their first contact with the
age barrier. Hence they yield a common admissible feedback.

Choose a short $h_1>0$ and use the local
intensity bound to make all prefixes $h\le h_1$ negative when $x$
is large. For $h\ge h_1$, if $F_\infty'(x)-1\le\delta$, then
\[
I_\infty(x,a;h)
\le\delta-q\int_0^{h_1}w_a(s)\,ds<0
\]
for sufficiently small $\delta$. This proves finiteness.
The same strict-prefix argument as in Theorem~\ref{thm:positive},
with \eqref{eq:inftytail} for long prefixes, proves continuity in age.
The infinite-horizon identity in Proposition~\ref{prop:dpp} is
attained by this explicit plan. Iterate it along the resulting barrier
strategy through $n$ claims. The unresolved continuation is at most
\[
e^{-qT_n}(x+cT_n+c/q)\mathbf1_{\{T_n<\tau^{L^b}\}}.
\]
It tends to zero almost surely and is bounded by a deterministic
constant, since $T_n\to\infty$ and $q>0$. Dominated convergence gives
$J^\infty(x,a;L^b)=V(x,a)$.
An infinite maximizing payment time means no active payment of that
layer on the no-claim branch, not a positive probability of no future
claim. This proves the proposition.
\end{proof}

\section{Nonexponential claims: failure of barrier optimality}
\label{sec:negative}
We prove Theorem~\ref{thm:conflict} by comparing two initial surplus
levels at the same time and claim age. The argument first
forces full payment at small surplus and then shows a strict benefit
from retention at a larger surplus. These two decisions exclude any
single upper barrier. Concavity and differentiability of the value
are not assumed.

\subsection{Properties of the claim-size distribution}
Define
\begin{equation*}
K(x)=\mu S(x)-\int_x^\infty S(y)\,dy.
\end{equation*}
Equivalently, $K(x)=S(x)\{\mu-\E[Y-x\mid Y>x]\}$.
It compares the mean of a fresh claim with the expected amount by
which a claim exceeds $x$, conditional on exceeding it. The factor
$S(x)$ weights this difference by the chance of such an exceedance.
The new-better-than-used-in-expectation (NBUE) property means
$\E[Y-x\mid Y>x]\le\mu$, or $K(x)\ge0$.

For the converse we need strictness of this comparison, not merely
its nonnegative sign. The next lemma also controls the moments of
remaining claim sizes and waiting times used in the subsequent bounds.
\begin{lemma}\label{lem:nbue}
For nonexponential $G\in\mathcal I$,
\begin{equation}\label{eq:strictnbue}
K(x)>0\quad(x>0),\qquad \mu g(0)<1.
\end{equation}
If a positive random variable $Z$ has mean residual life at most $L$,
then, for $n\ge1$ and $z\ge0$ with $\mathbb P(Z>z)>0$,
\begin{equation}\label{eq:moments}
\begin{aligned}
\E[(Z-z)^n\mid Z>z]&\le n!L^n,\\
\E[e^{r(Z-z)}\mid Z>z]&\le(1-Lr)^{-1},\qquad 0<r<L^{-1}.
\end{aligned}
\end{equation}
\end{lemma}
\begin{proof}
IFR gives $S(x+s)/S(x)\le S(s)$. If the two integrals agree for
some $x>0$, continuity forces equality of these survivals for every
$s$, hence $h_Y(x+s)=h_Y(s)$. A nondecreasing periodic function is
constant, contradicting nonexponentiality. This proves strict
positivity of $K$. If $h_Y(0)>0$, nonconstancy gives
$\mu<1/h_Y(0)$ by comparison with $e^{-h_Y(0)s}$; if $h_Y(0)=0$
the second inequality is immediate.

For $n\ge2$, Tonelli's theorem and the residual-mean bound give
\[
\begin{aligned}
\int_z^\infty(u-z)^{n-1}S_Z(u)\,du
&=(n-1)\int_z^\infty(v-z)^{n-2}
       \int_v^\infty S_Z(u)\,du\,dv\\
&\le(n-1)L\int_z^\infty(v-z)^{n-2}S_Z(v)\,dv .
\end{aligned}
\]
With $m_n(z)=\E[(Z-z)^n\mid Z>z]$, this is
$m_n(z)\le nL m_{n-1}(z)$. Since $m_1(z)\le L$,
\[
m_n(z)\le n!L^n,\qquad
\E[e^{r(Z-z)}\mid Z>z]
=\sum_{n=0}^\infty\frac{r^n m_n(z)}{n!}
\le\sum_{n=0}^\infty(rL)^n=\frac1{1-rL},
\]
which proves \eqref{eq:moments}.
\end{proof}

\subsection{Dividend bounds at small safety loadings}
Put
\begin{equation}\label{eq:scaling}
c_0=\mu/m,\qquad c=c_0+\varepsilon\sqrt q,\qquad 0<\varepsilon\le1.
\end{equation}
Constants $C$ below depend on the two fixed distributions but not on
age or sufficiently small $q$ and $\varepsilon\le1$.

We choose premium income only slightly above expected claims per
unit time. The next bound controls how many claims can contribute
to the dividend comparison before ruin.
\begin{lemma}\label{lem:count}
For $0<\varepsilon\le1$ and sufficiently small $q>0$, the
no-dividend process started at time zero has ruin time $\tau^0$
satisfying
\begin{equation}\label{eq:count}
C_q(x,a):=\E_{x,a}\sum_{T_n\le\tau^0}e^{-qT_n}
\le C(1+x)q^{-1/2},
\end{equation}
where the ruin-causing claim is included. Thus $C_q(x,a)$ is an
expected claim count with later claims receiving smaller weights;
it is a count, rather than a probability of ruin.
\end{lemma}
The proof is given in Appendix~\ref{app:count}. The estimate is
uniform in claim age, which is needed when the following dividend
bounds are applied after different waiting periods.

The count bound gives the following estimates for dividends and for
the loss associated with reserve held until the next claim. The second
estimate is the key comparison at small surplus.
\begin{lemma}[Dividend and retention bounds]\label{lem:energy}
Under \eqref{eq:strictnbue} and \eqref{eq:scaling},
\begin{equation}\label{eq:linearbound}
V(x,a)\le x+c m_A(a)+C\varepsilon(1+x).
\end{equation}
For a policy satisfying $J^T(t,x,a;L)\ge U(t,x,a)-\xi$, let $Z$ be
its surplus immediately before the first claim. Then
\begin{equation}\label{eq:energy}
\E[e^{-qA_a}K(Z);A_a\le T-t]
\le C\varepsilon(1+x)+Cq+cM e^{-q(T-t)}+\xi.
\end{equation}
The expectation measures the discounted loss $K(Z)$ at the first
claim; paths on which maturity comes first contribute zero.
On the infinite horizon, omit the last tail term and the event
$A_a\le T-t$.
\end{lemma}
\begin{proof}
The function $L_c(x,a)=x+c m_A(a)$, killed at ruin, satisfies
$m_A'=\eta m_A-1$ and
\begin{equation*}
\mathcal L L_c(x,a)
=\eta(a)\{(cm-\mu)G(x)-K(x)\},
\end{equation*}

Indeed,
\[
\begin{aligned}
\mathcal L L_c
&=c+c(\eta m_A-1)
 +\eta(a)\{\J_g L_c(x,0)-x-cm_A(a)\}\\
&=\eta(a)\left\{cmG(x)-\int_0^xS(y)\,dy\right\},\\
\int_0^x S(y)\,dy&=\mu G(x)+K(x).
\end{aligned}
\]
For a dividend $d\le x$,
$L_c(x-d,a)-L_c(x,a)=-d$. Write $\tau=\tau^L$ and $X=X^L$.
For $s\ge t$, with $s\le T$ in the finite-horizon problem, put
\[
J_L^s=\E_{t,x,a}\int_{[t,s\wedge\tau]}e^{-q(r-t)}\,dL_r.
\]
Applying the discounted finite-variation chain rule gives
\begin{align}
J_L^s
&+\E[e^{-q(s-t)}L_c(X_s,a_s);s<\tau]\nonumber\\
&+\E\sum_{T_n\le s\wedge\tau}e^{-q(T_n-t)}K(Z_n)
 +q\E\int_t^{s\wedge\tau}e^{-q(r-t)}L_c(X_r,a_r)\,dr\nonumber\\
&=L_c(x,a)
 +(cm-\mu)\E\sum_{T_n\le s\wedge\tau}e^{-q(T_n-t)}G(Z_n).
\label{eq:account}
\end{align}
Here expectations are under $\E_{t,x,a}$, $J_L^s$ is the expected
discounted dividends from $t$ to $s$, and
$Z_n$ is pre-claim surplus, including at ruin. The claim compensator
is $\eta(a_{s-})g(y)\,ds\,dy$ in the observed filtration.
At a claim, settlement precedes any surviving post-claim dividend;
the latter is accounted for by the exact linear payment identity.
Initial lump-sum payments and terminal liquidation have the same identity.

Any dividend policy ruins no later than the no-dividend process on
the same sample path. By time homogeneity, Lemma~\ref{lem:count}
bounds the positive term
in \eqref{eq:account} by $C\varepsilon(1+x)$.
All terms on the left except dividends are nonnegative.
Localizing by time, age and claim count and using $0\le K\le\mu$,
$m_A\le M$ and \eqref{eq:count} justifies expectations and limits.
The discounted terminal term vanishes on the infinite horizon because
it is bounded by $e^{-q(s-t)}(x+c(s-t)+cM)$.
This proves \eqref{eq:linearbound}.

For a nearly optimal policy, compare its payoff with
$x+c\int_0^{T-t}w_a(s)\,ds$ from \eqref{eq:budget} and retain only the first
$K$ term. The initial capital $x$ cancels. Finally,
\[
c m_A(a)-c\int_0^{T-t}w_a(s)\,ds
\le\tfrac12cq\,\E A_a^2+cM e^{-q(T-t)}
\le cM^2q+cM e^{-q(T-t)},
\]
by \eqref{eq:moments}. This proves \eqref{eq:energy}.
\end{proof}

\subsection{Optimal dividends at small surplus}
For $0\le r\le T$ set $F_r(z)=\J_g U(r,z,0)$.

The loss estimate is global. We next need a local comparison: at
small reserve, the increase in expected post-claim continuation is
strictly less than the extra capital retained. This will make certain
delayed dividends strictly improvable.
\begin{lemma}\label{lem:local-source}
There exist $r_0,\kappa>0$ and parameter bounds such that
\begin{equation}\label{eq:localslope}
F_r(v)-F_r(u)\le(1-\kappa)(v-u),
\qquad 0\le u<v\le r_0,\quad 0\le r\le T.
\end{equation}
The same bound holds for $\J_gV$. These constants do not depend on
the horizon, age or subsequently smaller parameters.
\end{lemma}
\begin{proof}
On a bounded interval, $h_Y$ is a continuous function of bounded
variation and $g=h_YS$ satisfies the measure identity
\[
dg=S\,dh_Y-h_Y^2S\,dy.
\]
For nonnegative $f$ bounded by $M_0$ on $[0,r_0]$, the convolution
has weak derivative
\begin{align}
(\J_g f)'(z)
&=g(0)f(z)+\int_{(0,z]}f(z-y)\,dg(y)\nonumber\\
&\le M_0\left[h_Y(0)+\int_{(0,z]}S(y)\,dh_Y(y)\right]
\le M_0h_Y(r_0).
\label{eq:BVbound}
\end{align}
This uses only a derivative measure of the density, not a derivative
of the value. Integrating the bound is legitimate since the convolution
is locally absolutely continuous. By \eqref{eq:linearbound},
\[
\sup_{r,\ z\le r_0}U(r,z,0)
\le r_0+\mu+m\varepsilon\sqrt q+C\varepsilon(1+r_0).
\]
The strict inequality $\mu h_Y(0)<1$ lets us choose $r_0,\kappa$
so that $h_Y(r_0)(\mu+2r_0)<1-2\kappa$, and then choose
$\varepsilon,q$ small enough that the displayed upper bound does
not exceed $\mu+2r_0$. Equations \eqref{eq:BVbound} and
\eqref{eq:linearbound} prove the assertion, also for $V$.
\end{proof}

The local bound allows us to compare two plans that differ only in
when a small dividend is paid, and then identify the optimal initial
payment.
From now on write $H(s,z)$ for the true continuation source in
\eqref{eq:H}. Equation \eqref{eq:localslope} implies
\begin{equation}\label{eq:Hdecrease}
H(s,v)-H(s,u)
\le-[q+\kappa\eta(a+s)](v-u)
\quad(0\le u<v\le r_0).
\end{equation}

\begin{lemma}[Advancing a payment]\label{lem:advance}
For positive remaining time, if an optimal claim-free plan retains
$r\in(0,r_0)$ immediately after its initial payment, it makes no
subsequent payment before first reaching $r_0$. On a finite horizon
it cannot remain below $r_0$ until maturity. The no-payment assertion
also holds on the infinite horizon.
\end{lemma}
\begin{proof}
For the original plan $\ell$, put $p_s=\ell_s-\ell_0$ and
$r=x-\ell_0$, so $p_0=0$ and $z_s=r+cs-p_s$.
Suppose $p_{s_1}>0$ before the first upward passage to $r_0$.
Choose $0<d<\min(r,p_{s_1})$ and advance this amount to the initial
time $t$, deducting it from the earliest subsequent payments:
\[
\widetilde\ell_s=\ell_0+\max(d,p_s),\qquad
\widetilde z_s=z_s-(d-p_s)^+.
\]
When $p_s<d$, the new surplus is $r-d+cs\ge0$; otherwise the two
paths have coalesced. Before coalescence both lie in $[0,r_0]$.
The first-claim objectives differ by at least
\begin{equation*}
\int_0^{s_1}w_a(s)[q+\kappa\eta(a+s)](d-p_s)^+\,ds>0.
\end{equation*}
Strictness follows from $q>0$ and right continuity of $p$ at zero.
This contradicts optimality. If the path stays below $r_0$ until
maturity, advance part of its terminal payment by the same construction.
The infinite-horizon assertion uses a finite coalescence time and has
no additional terminal issue.
\end{proof}

The preceding comparison constrains any maximizing plan that retains
positive capital. Combining it with the first-claim loss bound excludes
such retention and gives an exact identity for the value near zero.
\begin{proposition}\label{prop:small}
For sufficiently small fixed $\varepsilon$, then sufficiently small
fixed $q$, there are a compact age interval $J$, $r_0>0$ and
$H_*<\infty$ such that every optimal policy from
$a\in J$, $0<x\le r_0/2$, $T-t\ge H_*$ initially pays exactly $x$.
In particular,
\begin{equation}\label{eq:contact}
U(t,x,a)=x+U(t,0,a)
\quad(0\le x\le r_0/2,\ a\in J,\ T-t\ge H_*).
\end{equation}
The same assertions hold for $V$.
\end{proposition}
\begin{proof}
A proper waiting-time law has a continuous hazard that is not identically
zero. Choose an age interval $J$ with a slightly larger neighborhood
on which $0<\eta_*\le\eta\le\eta^*$. Put
$c_{\min}=c_0$, $c_{\max}=c_0+1$ and, after decreasing $r_0$ if
necessary, $s_0=r_0/(4c_{\max})$, so the hazard bounds hold on
$a+[0,s_0]$ for every $a\in J$.
The constants
\[
k_*=\min_{c_{\min}s_0/2\le z\le r_0}K(z)>0,\qquad
d_*=\tfrac12s_0\eta_*e^{-(q_0+\eta^*)s_0}k_*>0
\]
are now fixed independently of the later choices of $\varepsilon,q$.

Fix a state and take a maximizing deterministic first-claim plan,
which exists by Proposition~\ref{prop:dpp}. By
\eqref{eq:near-continuation}, it can be extended to a $\delta$-optimal
policy without changing its first-claim reserve path. Lemma~\ref{lem:energy}
applies to these extensions. Letting $\delta\downarrow0$ gives
\[
\int_0^{T-t}w_a(s)\eta(a+s)K(z_s)\,ds
\le C\varepsilon(1+x)+Cq+cMe^{-q(T-t)}.
\]
If the maximizing plan retains $r\in(0,x]$ initially, Lemma~\ref{lem:advance}
makes its no-claim capital $r+cs$ until it reaches $r_0$.
For $s\le s_0$ this is at most $3r_0/4$; for $s\ge s_0/2$ it is
at least $c_{\min}s_0/2$. Hence, with $T-t\ge r_0/c_{\min}$,
\begin{equation}\label{eq:loss-lower}
\E[e^{-qA_a}K(Z);A_a\le T-t]
\ge\int_{s_0/2}^{s_0}w_a(s)\eta(a+s)K(r+cs)\,ds\ge d_*.
\end{equation}
Choose $\varepsilon$ so $C\varepsilon(1+r_0/2)<d_*/3$,
then $q$ so $Cq<d_*/3$, and finally $H_*$ so
$cM e^{-qH_*}<d_*/3$. Also impose all previous parameter bounds
and $H_*\ge r_0/c_{\min}$. Lemma~\ref{lem:energy} contradicts
\eqref{eq:loss-lower}. Thus every maximizing plan pays $x$ initially. Subtracting this
payment leaves a feasible plan from zero, so \eqref{eq:firstclaim}
gives $U(t,x,a)\le x+U(t,0,a)$. The reverse inequality follows by
immediate payment, proving \eqref{eq:contact}.

For any optimal admissible policy $L$, with initial no-claim plan
$\ell$, conditioning gives
\[
U(t,x,a)=J^T(t,x,a;L)\le\Phi_U(\ell)\le U(t,x,a).
\]
Its plan is therefore maximizing and $L_t=x$. On the infinite horizon,
\eqref{eq:tails} passes the contact identity to $V$. If an infinite-horizon
optimal policy retained positive capital, its first-claim plan would
again be maximizing; Lemmas~\ref{lem:advance} and \ref{lem:energy},
with no terminal tail, give the same contradiction.
\end{proof}

All constants $r_0,\kappa,J,s_0,d_*$ can be fixed from the two laws,
using a preliminary upper bound $q\le q_0$. Then $\varepsilon$ is
fixed and $q$ is chosen sufficiently small. If $q$ is decreased in
the comparison below, the horizon cutoff is increased again; no
constant in the local loss lower bound is changed.

\subsection{The benefit of retaining a larger surplus}
It remains to show that full distribution fails at a larger surplus.
We first identify the value of continuing from zero over a short
interval. This provides the common comparison term for the two initial
decisions.

\begin{lemma}\label{lem:zero-prefix}
Choose $a_*$ in the interior of $J$ and $h_0>0$ such that
$a_*+[0,h_0]\subset J$ and $c_{\max}h_0\le r_0/2$.
For $B(t,a)=U(t,0,a)$ and $T-t\ge H_*+h_0$,
\begin{equation}\label{eq:zero-prefix}
B(t,a_*)=
c\int_0^{h_0}w_{a_*}(s)\,ds+
w_{a_*}(h_0)B(t+h_0,a_*+h_0).
\end{equation}
Thus keeping zero surplus and paying premium income is optimal on
this short interval from zero. The analogous identity holds for $V(0,a)$.
\end{lemma}
\begin{proof}
Every no-claim path from zero stays below $cs\le r_0/2$ before $h_0$.
The terminal continuation at $h_0$ equals
$z+B(t+h_0,a_*+h_0)$ by \eqref{eq:contact}.

The DPP stopped at the first claim or $h_0$ gives
\[
\begin{aligned}
B(t,a_*)=\sup_\ell\bigg\{
&\int_{[0,h_0]}w_{a_*}(s)\,d\ell_s
 +\int_0^{h_0}w_{a_*}(s)\eta(a_*+s)F_{t+s}(z_s)\,ds\\
&+w_{a_*}(h_0)\big[z_{h_0}+B(t+h_0,a_*+h_0)\big]\bigg\}.
\end{aligned}
\]
Integration by parts cancels the terminal term
$w_{a_*}(h_0)z_{h_0}$, leaving
\[
\begin{aligned}
B(t,a_*)&=c\int_0^{h_0}w_{a_*}(s)\,ds
+w_{a_*}(h_0)B(t+h_0,a_*+h_0)\\
&\quad+\sup_z\int_0^{h_0}w_{a_*}(s)H(s,z_s)\,ds .
\end{aligned}
\]
By \eqref{eq:Hdecrease},
$H(s,z)\le-[q+\kappa\eta(a_*+s)]z\le0$ on the relevant range.
The zero path attains the supremum, proving \eqref{eq:zero-prefix}.
\end{proof}

\begin{proof}[Proof of Theorem~\ref{thm:conflict}]
Fix the small $\varepsilon$ selected above, and reduce $q$ further
as necessary. The feasible zero-barrier strategy gives
\[
U(r,z,0)\ge z+b_0(r),\qquad
b_0(r)=c\int_0^{T-r} e^{-qs}S_0(s)\,ds.
\]

Since $\int_0^z(z-y)g(y)\,dy=z-\int_0^zS(y)\,dy$,
\begin{equation*}
\begin{aligned}
F_r(z)-z
&\ge b_0(r)G(z)-\int_0^zS(y)\,dy\\
&=(b_0(r)-\mu)G(z)-K(z).
\end{aligned}
\end{equation*}
Moreover
\[
b_0(r)-\mu
\ge m\varepsilon\sqrt q-Cq-cM e^{-q(T-r)}.
\]
Increase $H_*$ and reduce $q$ so this is at least
$m\varepsilon\sqrt q/2$ for $T-r\ge H_*$.
These changes preserve Proposition~\ref{prop:small}.
At the final choice of parameters take the cutoff to satisfy both
$cMe^{-qH_*}<d_*/3$ and
$cMe^{-qH_*}<m\varepsilon\sqrt q/4$, as well as all earlier
lower bounds on $H_*$. The premium--discount pair is fixed only
after these simultaneous requirements and the inequalities below
have been imposed.

Set $\ell_*=\eta_*m\varepsilon/16$ and $x_q=\ell_*/\sqrt q$.
Since $S$ is nonincreasing and integrable, $xS(x)\to0$.
As $0\le K(x)\le\mu S(x)$, uniformly for $0\le s\le h_0$,
\[
G(x_q+cs)\to1,\quad K(x_q+cs)=o(\sqrt q),\quad
q(x_q+cs)=\ell_*\sqrt q+O(q).
\]
Thus $q$ can be fixed small enough that, whenever
$T-t\ge H_*+h_0$,
\begin{equation*}
\eta(a_*+s)\{F_{t+s}(x_q+cs)-(x_q+cs)\}
-q(x_q+cs)\ge \tfrac18\eta_*m\varepsilon\sqrt q>0
\end{equation*}
on $[0,h_0]$. For example use $G\ge3/4$,
$\eta^*K\le\eta_*m\varepsilon\sqrt q/8$ and
$q(x_q+cs)\le\eta_*m\varepsilon\sqrt q/8$.

The supremum obtainable after immediate full distribution is
$x_q+B(t,a_*)$. Instead, wait without dividends until the first claim
or $h_0$; on the no-claim branch at $h_0$, distribute $x_q+ch_0$.
Using the value function for subsequent decisions gives the comparison
quantity
\[
\begin{aligned}
W_{\rm wait}={}&\int_0^{h_0}
w_{a_*}(s)\eta(a_*+s)F_{t+s}(x_q+cs)\,ds\\
&+w_{a_*}(h_0)\{x_q+ch_0+B(t+h_0,a_*+h_0)\}.
\end{aligned}
\]
By \eqref{eq:near-continuation} and its stopped version, this quantity
is achievable within any error $\delta>0$. Hence
$U(t,x_q,a_*)\ge W_{\rm wait}$. Using \eqref{eq:zero-prefix} and
\[
w_{a_*}(h_0)(x_q+ch_0)-x_q
=\int_0^{h_0}w_{a_*}(s)
[c-(q+\eta(a_*+s))(x_q+cs)]\,ds,
\]
we find that $W_{\rm wait}-x_q-B(t,a_*)$ equals
\begin{equation}\label{eq:strictgain}
\int_0^{h_0}w_{a_*}(s)
\big[\eta(a_*+s)\{F_{t+s}(x_q+cs)-(x_q+cs)\}
-q(x_q+cs)\big]\,ds>0.
\end{equation}
This proves $U(t,x_q,a_*)>x_q+U(t,0,a_*)$.
Choose $x_*=r_0/2$, $x^*=x_q>x_*$ and $H_0=H_*+h_0$.
More precisely, the gain in \eqref{eq:strictgain} is at least
\[
\gamma:=\frac18\eta_*m\varepsilon\sqrt q
\int_0^{h_0}w_{a_*}(s)\,ds>0,
\]
independently of all sufficiently long maturities. Therefore
\eqref{eq:tails} passes both the contact identity and this strict gain
to $V$, at the same parameters and surplus levels.

Put $R(x)=U(t,x,a_*)-x$. Initial payment implies that $R$ is
nondecreasing. Equations \eqref{eq:contact} and \eqref{eq:strictgain}
give
\[
\frac{R(x_*)-R(0)}{x_*}=0
<
\frac{R(x^*)-R(0)}{x^*},
\]
which contradicts concavity. The same calculation applies to $V-x$.

Choose $t_*\ge a_*$ and $T_0=t_*+H_0$. For every $T\ge T_0$,
$T-t_*\ge H_0$, so \eqref{eq:main-small}--\eqref{eq:main-large}
hold at physical states. If $b_T(t_*,a_*)>0$, then for
$0<x<\min\{x_*,b_T(t_*,a_*)\}$ its initial payment is
\[
(x-b_T(t_*,a_*))^+=0<x,
\]
contradicting \eqref{eq:main-small}. If $b_T(t_*,a_*)=0$, its payoff
from $x^*$ is at most
\[
x^*+U^T(t_*,0,a_*)<U^T(t_*,x^*,a_*),
\]
contradicting optimality. The argument for an age barrier is identical.
This proves (i)$\Rightarrow$(ii) and the stronger assertions.
Conversely, if $G$ were exponential,
Theorem~\ref{thm:positive} would give a globally optimal barrier for
every $c,q,T>0$, contradicting (ii).
\end{proof}

\begin{proof}[Proof of Theorem~\ref{thm:classification}]
Theorem~\ref{thm:positive} proves (i)$\Rightarrow$(ii), and
Proposition~\ref{prop:infinite-positive} proves
(i)$\Rightarrow$(iii). If (i) fails,
Theorem~\ref{thm:conflict} supplies one pair $(c,q)$ excluding
every globally optimal infinite-horizon barrier and every globally
optimal finite-horizon barrier for $T\ge T_0$.
Thus neither (ii) nor (iii) can hold, proving both reverse implications.
\end{proof}

\section{Claim distributions and numerical examples}\label{sec:examples}
\subsection{Examples and scope}
Gamma claims with shape $\alpha>1$ and Weibull claims with shape
$k>1$ belong to $\mathcal I$. Theorem~\ref{thm:conflict} therefore
applies to both families for every fixed arrival law satisfying
Assumption~\ref{ass:arrival}. Shape one gives the exponential case
of Theorem~\ref{thm:classification}. The Weibull hazard is unbounded;
the local estimate \eqref{eq:BVbound} is sufficient for this example.

The IFR restriction matters. Nontrivial exponential mixtures fall
under established Poisson/L\'evy sufficient conditions for barrier
optimality \citep{Loeffen2008,LoeffenRenaud2010}, and lie outside
the continuous IFR class. Within that class, nonexponential claims
may still admit an optimal barrier at particular parameters: if the
claim hazard is bounded by $H$, the zero barrier is optimal whenever
$q\ge cH$. A verification and an extension of the converse beyond
IFR are given in \suppref{supp:zero} and \suppref{supp:extension}.

\subsection{Dividend decisions for exponential and Gamma claims}\label{sec:numerics}
Take independent interclaim times and claim sizes with
$A\sim\operatorname{Gamma}(2,2)$ and $Y\sim\operatorname{Gamma}(2,2)$,
where the second argument is the rate. Both means are one. Set
$c=1.05$, $q=0.01$, $T=30.5$, and examine the state $t=a=0.5$.
For comparison, replace only the claim-size law by $\Exp(1)$.
Time and money are measured in units of mean interclaim time and
mean claim size. Thus the safety loading is $5\%$ in both cases.

With time and age step $h$, surplus step $\delta$, and
$t_k=T-kh$, write $u_k(x,a)\approx U(t_k,x,a)$, with $u_0(x,a)=x$.
The scheme first computes the value of waiting one time step and
then optimizes the initial payment over all surplus grid points:
\begin{align}
Q_k(x,a)={}&w_a(h)u_{k-1}(x+ch,a+h)\nonumber\\
&+\int_0^h w_a(s)\eta(a+s)
       \J_g u_{k-1}(x+cs,0)\,ds,\label{eq:numerical-Q}\\
u_k(x_i,a)={}&x_i+\max_{0\le j\le i}
                   \{Q_k(x_j,a)-x_j\}.\label{eq:numerical-envelope}
\end{align}
Here $\J_g u_{k-1}(z,0)=\int_0^z u_{k-1}(z-y,0)g(y)\,dy$.
In the claim term, $t_k+s$ is approximated by $t_{k-1}$; hence
\eqref{eq:numerical-Q} is a time discretization. No barrier form is
imposed in \eqref{eq:numerical-envelope}. Quadrature, interpolation
and boundary treatment are specified in \suppref{supp:numerics}.

Let $y_{h,\delta}(x)$ be the selected post-payment reserve at
$t=a=0.5$, choosing the smallest reserve in a tie, and set
\[
\Delta_{h,\delta}(x)=u_{h,\delta}(0.5,x,0.5)
-x-u_{h,\delta}(0.5,0,0.5),
\]
where $u_{h,\delta}$ is the time interpolant of the grid values.
This is the computed gain over immediate full payment.
Table~\ref{tab:numerics} reports the end $x_0$ of the initial grid
interval with $y_{h,\delta}=0$, the largest selected reserve
$r_{\max}$, and $\Delta_{h,\delta}(3)$.

\begin{table}[tb]
\centering\small
\caption{Refinement at $T=30.5$, $t=a=0.5$, $c=1.05$, $q=0.01$,
and $x_{\max}=12$. Both claim laws have mean one.}
\label{tab:numerics}
\begin{tabular}{llrrrr}
\toprule
Claims & $h=\delta$ & $u_{h,\delta}(t,0,a)$ & $x_0$ & $r_{\max}$ & $\Delta(3)$ \\
\midrule
Gamma(2,2) & 0.0500 & 0.75876 & 1.5000 & 3.6500 & 0.03291 \\
Gamma(2,2) & 0.0250 & 0.76986 & 1.5750 & 3.6500 & 0.03135 \\
Gamma(2,2) & 0.0125 & 0.77590 & 1.5875 & 3.6375 & 0.03049 \\
\addlinespace
Exponential(1) & 0.0500 & 0.83922 & 0.0000 & 3.1000 & 0.07176 \\
Exponential(1) & 0.0250 & 0.83937 & 0.0000 & 3.1250 & 0.07194 \\
Exponential(1) & 0.0125 & 0.83943 & 0.0000 & 3.1250 & 0.07203 \\
\bottomrule
\end{tabular}

\end{table}

\begin{figure}[tb]
\centering
\includegraphics[width=\textwidth]{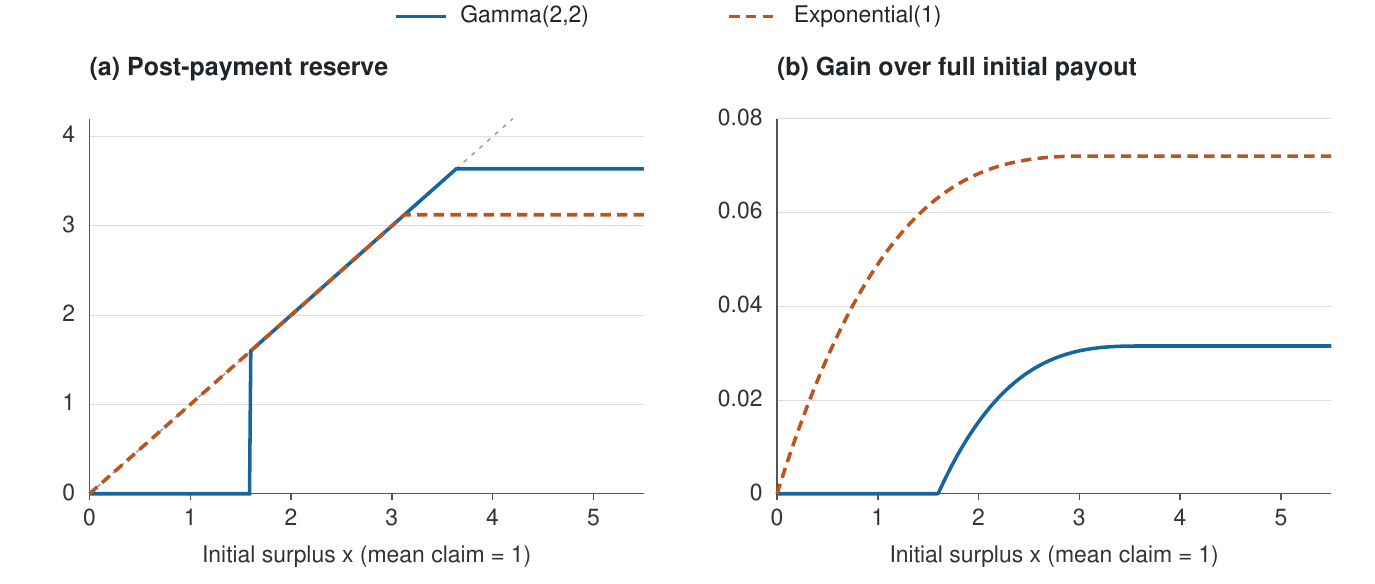}
\caption{Computed post-payment reserves (left) and gains relative to
full initial distribution (right), with $h=\delta=0.0125$.
The diagonal in the left panel denotes no initial payment.}
\label{fig:numerics}
\end{figure}

For Gamma claims, Figure~\ref{fig:numerics} shows full payment at
$x=0.5$ and retention of the entire surplus at $x=3$. On the finest
grid, $x_0\simeq1.59$ and $\Delta_{h,\delta}(3)\simeq0.0305$.
These decisions cannot be produced by one upper barrier: full
payment at a positive surplus requires barrier zero, which would
also pay the larger surplus in full. For exponential claims,
the selected reserve is capped near $3.13$, with no positive initial
full-payment interval, consistent with Theorem~\ref{thm:positive}.

The refinement in Table~\ref{tab:numerics}, separate changes of the
time and surplus steps, and an enlarged surplus boundary preserve
these decisions; the additional checks are reported in
\suppref{supp:numerics}. The computations illustrate how changing the
claim-size distribution at fixed mean can change the form of the
payment rule. They provide numerical evidence at the stated
parameters, rather than an error bound for the continuous problem.
A separate Gamma example with fully explicit analytical bounds is
given in \suppref{app:certificate}.

\Needspace{7\baselineskip}
\section{Conclusion}\label{sec:conclusion}
The claim-size distribution can determine whether a simple dividend
barrier is sufficient for optimal capital distribution. For a fixed
renewal law satisfying the stated assumptions, exponential claims are
exactly the distributions within the continuous IFR class that admit
an optimal barrier for every positive safety loading, discount rate
and finite horizon. Every nonexponential law in this class instead
admits one fixed premium--discount pair at which barrier optimality
fails for all sufficiently long horizons. The same distinction holds
on the infinite horizon.

The failure is explained by opposing initial decisions at the same
time and claim age: full payment at small surplus and a strict benefit
from retention at a larger surplus. The numerical comparison shows
this pattern for Gamma claims at moderate parameter values, while
exponential claims with the same mean give the expected barrier form.
This identifies a structural consequence of choosing a claim-size law
that cannot be seen from its mean alone.

The bounded residual waiting-time assumption implies the existence
of a positive exponential moment and excludes heavy-tailed arrivals
without such moments. Extending the characterization
to heavier-tailed renewal arrivals, determining the complete optimal
strategy when a barrier fails, and adding explicit solvency constraints
are natural further questions. Quantitative error bounds for the
numerical scheme would also strengthen the comparison of dividend
rules at fixed parameters.

\section*{Funding}
This work was supported in part by the National Natural Science Foundation
of China (NSFC) under Grant Nos. 12271274 and 12571514.

\section*{Declaration of competing interest}
The author declares no known competing financial interests or personal
relationships that could have appeared to influence the work reported
in this paper.

\appendix
\renewcommand{\thetheorem}{\Alph{section}.\arabic{theorem}}
\renewcommand{\thelemma}{\thetheorem}
\renewcommand{\theproposition}{\thetheorem}
\section{Payment plans before the next claim}
\label{app:plans}
This appendix supplies the continuity and finite-horizon maximization
facts used in Proposition~\ref{prop:dpp}. Its measurable choices are
used for truncated problems and near-optimal continuation.
\begin{lemma}\label{lem:plans}
For a continuous continuation source, the finite-horizon first-claim
objective has an attained maximum, a continuous value for $t<T$,
and a Borel maximizing plan. Terminal liquidation is specified
separately at $t=T$.
\end{lemma}
\begin{proof}
Fix $t<T$, write $d=x+c(T-t)$ and $\alpha=x/d\in[0,1)$,
and identify a plan with a probability measure $\nu$ on $[0,1]$ by
$\ell_{(T-t)v}=d\nu([0,v])$. The feasible correspondence is
\begin{equation*}
\mathcal K_\alpha=
\{\nu\in\mathcal P([0,1]):
\nu([0,v])\le\alpha+(1-\alpha)v\text{ for all }v\in[0,1]\}.
\end{equation*}
The space $\mathcal P([0,1])$ with weak convergence is compact
metrizable, and $\mathcal K_\alpha$ is nonempty (it contains
$\delta_1$). To check its closed graph, take
$\nu_j\Rightarrow\nu$, $\alpha_j\to\alpha$. For any $v<1$,
apply the bound at a continuity point $u>v$ of $\nu$, pass to the
limit, and then let $u\downarrow v$. This gives the constraint at $v$,
including $v=0$; the constraint at one is automatic.

The correspondence is also lower hemicontinuous. Given
$\nu\in\mathcal K_\alpha$ and $\alpha_j\to\alpha$, use $\nu_j=\nu$
if $\alpha_j\ge\alpha$. If $\alpha_j<\alpha$, push $\nu$ forward
by the map
\begin{equation*}
\phi_j(v)=\frac{\alpha-\alpha_j+(1-\alpha)v}{1-\alpha_j}.
\end{equation*}
For $u<\phi_j(0)$ its cumulative mass is zero; otherwise applying
the constraint at $\phi_j^{-1}(u)$ proves
$\nu_j([0,u])\le\alpha_j+(1-\alpha_j)u$.
Since $\phi_j\to{\rm id}$ uniformly, $\nu_j\Rightarrow\nu$.
This establishes continuity of the feasible correspondence.

After $s=(T-t)v$, the nonconstant objective is
\[
(T-t)\int_0^1 w_a((T-t)v)
H_f(t,a;(T-t)v,x+c(T-t)v-d\nu([0,v]))\,dv.
\]
Under simultaneous convergence of $(t,x,a,\nu)$, distribution
functions converge outside the countable atom set of the limiting
measure. The integrand converges almost everywhere and is bounded
on compact parameter sets. Dominated convergence therefore gives
joint continuity on the feasible graph. If a function on the full
product space is needed, replace its surplus argument by its positive
part; this is a continuous extension for the integrated objective.
The compact metrizable choice space, continuous feasible
correspondence and continuous objective meet the hypotheses of
Berge's theorem and the measurable maximum theorem
\citep[Theorem~18.19]{AliprantisBorder2006}.
The map $(\nu,v)\mapsto\nu([0,v])$ is Borel, so the selector gives
a jointly measurable right-continuous cumulative payment plan.

At $t=T$ select the initial payment $x$ separately.
For the recursively defined values in Proposition~\ref{prop:dpp},
the discounted capital budget bounds the difference from $x$ by
$c(T-t)$, proving continuity also at this boundary. Concatenation after
a claim is measurable because the residual state and the selected
cumulative plan are Borel; claim times have no finite accumulation.

\end{proof}

\section{Discounted claims before ruin}\label{app:count}
This appendix proves Lemma~\ref{lem:count}, which is used only through
its uniform bound in the dividend estimates of Section~\ref{sec:negative}.
\begin{proof}
For generic $A\sim F_A$ and $Y\sim G$, write
$\varphi(s)=\E e^{-sA}$, $\varphi_a(s)=\E e^{-sA_a}$ and
$M_Y(r)=\E e^{rY}$. The first two transforms are the expected
discount factors over a fresh and a residual waiting time, respectively;
the last records exponential moments of a claim.
There is a positive root $r=r_q\le C\sqrt q$ of
\begin{equation}\label{eq:root}
M_Y(r)\varphi(q+cr)=1.
\end{equation}
Indeed, \eqref{eq:moments} gives fixed neighborhoods of exponential
integrability for both laws. Taylor expansion yields
\begin{align}
\log M_Y(r)+\log\varphi(q+cr)
&=-mq-m\varepsilon\sqrt q\,r+\tfrac12\sigma_0^2r^2
+R(q,r,\varepsilon),\nonumber\\
\sigma_0^2&=\Var Y+c_0^2\Var A>0,\nonumber\\
|R(q,r,\varepsilon)|
&\le C(q^2+qr+r^3+\sqrt q\,r^2).\nonumber
\end{align}
At $r=C_1\sqrt q$ a sufficiently large fixed $C_1$ makes the
expression positive, uniformly in $\varepsilon\le1$, while at zero
it is negative. Choose a root and reduce $q$ so that $r<1/(2\mu)$.

Set $O=-X^0_{\tau^0}$ on $\{\tau^0<\infty\}$ and $O=0$ otherwise,
and put $\chi_q=\E e^{-q\tau^0}$, with the exponential set to zero
at infinity.
Thus $O$ is the deficit at ruin, and $\chi_q$ is the present value of
one unit paid at ruin, averaged over paths; survival forever contributes
zero.
Stopping the exponential claim-epoch martingale gives
\begin{equation}\label{eq:ruintransform}
\E[e^{-q\tau^0+rO};\tau^0<\infty]
=e^{-rx}M_Y(r)\E e^{-(q+cr)A_a}.
\end{equation}

To justify the limit, let $\nu=\inf\{n\ge1:X_{T_n}^0<0\}$.
The root identity implies
\[
\E\!\left[
e^{-q(T_{n+1}-T_n)-r(X_{T_{n+1}}^0-X_{T_n}^0)}
 \mid\mathcal F_{T_n}^{0,a}\right]
=M_Y(r)\varphi(q+cr)=1.
\]
Stopping at $\nu\wedge n$ yields
\[
\begin{aligned}
&\E[e^{-q\tau^0+rO};\nu\le n]
 +\E[e^{-qT_n-rX_{T_n}^0};\nu>n]\\
&\hspace{25mm}=e^{-rx}M_Y(r)\varphi_a(q+cr).
\end{aligned}
\]
The second term is at most
$\E e^{-qT_n}=\varphi_a(q)\varphi(q)^{n-1}\to0$.
Monotone convergence gives \eqref{eq:ruintransform}.
The residual moment bound gives
\[
e^{-rx}M_Y(r)\varphi_a(q+cr)
\le\frac{\chi_q}{1-\mu r},\qquad
\chi_q\ge(1-\mu r)e^{-rx}\varphi_a(q+cr).
\]
Using $1-e^{-rx}\le rx$ and
$1-\varphi_a(q+cr)\le(q+cr)m_A(a)$, we obtain
\[
1-\chi_q\le r(x+\mu)+(q+cr)M\le C(1+x)\sqrt q.
\]
The exact delayed count identity is
\[
C_q=\varphi_a(q)+\varphi(q)(C_q-\chi_q),\qquad
C_q=\frac{\varphi_a(q)-\varphi(q)\chi_q}{1-\varphi(q)}.
\]
The term $C_q-\chi_q$ counts the discounted claims survived, each
followed by a fresh waiting time. All sums are bounded by the
unrestricted discounted renewal count. Since $\varphi_a(q)\le1$,
\[
C_q\le 1+\frac{1-\chi_q}{1-\varphi(q)}
\le C(1+x)q^{-1/2},
\]
where the last step uses $1-\varphi(q)\sim mq$ and the preceding
bound on $1-\chi_q$.
\end{proof}

\section{Numerical implementation and additional checks}\label{supp:numerics}
This section gives the implementation details for the numerical
comparison in Section~\mainref{sec:numerics}. The model, parameters,
and definitions of $y_{h,\delta}$ and $\Delta_{h,\delta}$ are those
in that section. The comparison uses the same arrival law and mean
claim size for both claim distributions.

Positive three-point Gaussian quadrature is used in the time integral
of \mainequation{eq:numerical-Q}, and piecewise-linear product
integration is used in the claim convolution. Values between surplus
nodes are linearly interpolated. The convolution is checked on affine
functions against the identity
\[
\int_0^x(x-y+B)g(y)\,dy
=x-\int_0^x S(y)\,dy+B G(x),\qquad B\ge0.
\]
All integral and interpolation weights are nonnegative before
floating-point rounding. The maximization in
\mainequation{eq:numerical-envelope} allows every grid post-payment
reserve, including full and partial distribution; ties are resolved
by the smallest retained reserve.

The age grid follows $a\mapsto a+h$, starting with a sufficiently
wide terminal grid and shrinking at each time step. No artificial
age boundary is imposed. At the upper surplus boundary $x_{\max}$,
values beyond the grid are extended by
\[
u_k(x_{\max}+y,a)=u_k(x_{\max},a)+y,\qquad y\ge0.
\]
The main comparison uses $x_{\max}=12$; enlarging this boundary to
$18$ at $h=\delta=0.025$ changes the reported values by less than
$10^{-10}$.

Time and surplus refinement were also separated. For Gamma claims,
$(h,\delta)=(0.0125,0.025)$ and $(0.025,0.0125)$ give
$\Delta_{h,\delta}(3)=0.03010$ and $0.03140$, respectively, and
preserve full payment at $x=0.5$ and retention at $x=3$.
Together with the refinement table, these checks support numerical
stability at the stated parameters. They do not establish an error
bound for the continuous control problem. The analytical example in
Section~\ref{app:certificate} uses separate, conservative parameters.

\section{A zero barrier at large discount rates}\label{supp:zero}
Nonexponential claims can admit an optimal barrier at particular
parameters. The following verification explains why the quantifiers
in Theorem~\mainref{thm:classification} are essential.
\begin{proposition}\label{prop:zero-supp}
Under the model assumptions, suppose that the claim-size hazard
satisfies $h_Y(x)\le H<\infty$ for every $x\ge0$. If $q\ge cH$,
the zero barrier is optimal on both finite and infinite horizons.
\end{proposition}
\begin{proof}
On the finite horizon, its value is
\[
W(t,x,a)=x+c\int_0^{T-t}w_a(s)\,ds.
\]
The payment at time $t$ is $x$, and all premiums are paid until the
first claim or maturity. Direct differentiation gives $W_x=1$,
$W(T,x,a)=x$, and
\[
W_t+\mathcal LW-qW
=-qx+\eta(a)\left[cG(x)\int_0^{T-t}w_0(s)\,ds-\int_0^x S(y)\,dy\right]
\le0,
\]
because
\[
G(x)\le H\int_0^x S(y)\,dy,\qquad
c\int_0^{T-t}w_0(s)\,ds\le c/q,\qquad cH/q\le1.
\]
Extend $W$ by zero after ruin. The discounted finite-variation
formula, localized by claim count, then gives, for every admissible
$L$,
\[
J^T(t,x,a;L)+\E_{t,x,a}[e^{-q(T-t)}X_T^L;\,T<\tau^L]
\le W(t,x,a).
\]
Here dividends reduce $W$ by exactly their amount, and the
nonnegative terminal term can be discarded. The zero-barrier
payoff equals $W$, proving optimality. On the infinite horizon,
replace $T-t$ by $\infty$; the residual terminal term is at most
$e^{-qs}(x+cs+c/q)\to0$ when the initial time is zero.
\end{proof}

\section{An extension beyond increasing claim hazards}\label{supp:extension}
The converse proof uses two strict inequalities for the claim-size
law and local regularity of its density. They can hold without IFR.
\begin{corollary}\label{cor:NBUE}
Retain the arrival assumptions of Theorem~\mainref{thm:classification}.
The negative conclusions of Theorem~\mainref{thm:conflict} remain valid
for a continuous claim density $g$, with $S>0$ and $0<\mu<\infty$,
satisfying $K(x)>0$ for every $x>0$ and $\mu g(0)<1$, provided
$g$ is of bounded variation on a neighborhood of zero. Here
$K(x)=\mu S(x)-\int_x^\infty S(y)\,dy$.
\end{corollary}
\begin{proof}
The condition $K(x)>0$ implies that every residual claim has mean at
most $\mu$. Lemma~\mainref{lem:nbue} therefore supplies the moment
bounds used in Lemmas~\mainref{lem:count} and \mainref{lem:energy}.
For nonnegative $f$, the replacement for \mainequation{eq:BVbound} is
\[
(\J_gf)'(z)\le\sup_{[0,r_0]}f\,\{g(0)+|dg|((0,r_0])\},\qquad 0<z<r_0.
\]
Since $|dg|((0,r_0])\to0$ and $\mu g(0)<1$, choose $r_0,\kappa>0$
with
\[
(\mu+2r_0)\{g(0)+|dg|((0,r_0])\}<1-2\kappa.
\]
For the parameter choices in Lemma~\mainref{lem:local-source},
$\sup_{r,\,z\le r_0}U(r,z,0)\le\mu+2r_0$. Integrating the derivative
bound gives \mainequation{eq:localslope}. The payment comparison and
strict-gain argument of Theorem~\mainref{thm:conflict} then apply.
\end{proof}

\section{An explicit Gamma example}\label{app:certificate}
The following example makes the parameter choice in
Theorem~\mainref{thm:conflict} fully explicit. The bounds are conservative
and serve an analytical purpose; the numerical comparison in
Section~\mainref{sec:numerics} uses moderate parameters.
\begin{proposition}\label{prop:certificate}
Let both $A$ and $Y$ have Gamma shape $2$ and rate $2$, independently.
For
\begin{equation*}
q=10^{-24},\qquad c=1+10^{-20},
\end{equation*}
there is no globally optimal infinite-horizon age barrier and no
globally optimal finite-horizon time--age barrier for
$T\ge10^{26}+1$.
\end{proposition}
\begin{proof}
For the two Gamma laws,
\begin{gather*}
m=\mu=M=1,\qquad m_A(a)=\frac{a+1}{1+2a},\\
g(x)=4xe^{-2x},\qquad S(x)=(1+2x)e^{-2x},\qquad K(x)=xe^{-2x}.
\end{gather*}
Put $\varepsilon=10^{-8}$, so $c-1=\varepsilon\sqrt q=10^{-20}$.
The transforms are $M_Y(r)=(2/(2-r))^2$ and
$\varphi(s)=(2/(2+s))^2$. The positive root in \mainequation{eq:root}
satisfies
\[
cr^2-(2(c-1)-q)r-2q=0,\qquad r\le4\sqrt q.
\]
For the stated $q$, this lies inside the required transform domain.
The overshoot argument gives
$1-\chi_q\le(4x+13)\sqrt q$, whereas
$1-\varphi(q)\ge4q/9$. The delayed count identity implies
\[
C_q(x,a)\le32(1+x)/\sqrt q,\qquad
V(x,a)\le x+c m_A(a)+32\varepsilon(1+x).
\]
For a maximizing finite-horizon payment plan,
\begin{equation}\label{eq:certenergy}
\E[e^{-qA_a}K(Z);A_a\le T-t]
\le32\varepsilon(1+x)+2q+2e^{-q(T-t)}.
\end{equation}

Take $r_0=1/20$, $J=[1/2,1]$ and $s_0=1/160$.
The local continuation slope is bounded by
\[
\frac{4r_0}{1+2r_0}
\{r_0+c+32\varepsilon(1+r_0)\}<0.192,
\]
so $\kappa=1/2$ is valid. For $a\in J$ and $s\le s_0$,
$1\le\eta(a+s)\le2$. In \mainequation{eq:loss-lower} we can take
\[
k_*=\frac1{320}e^{-1/160},\qquad
d_*=\frac1{102400}e^{-1/40}>9.5\cdot10^{-6}.
\]
For $x\le1/40$, $T-t\ge10^{26}$, the upper bound
\eqref{eq:certenergy} is less than $3.3\cdot10^{-7}$.
Thus full initial payment is necessary for optimality on this small interval.

For the larger surplus choose $a_*=1/2$, $h_0=1/200$ and $x^*=30$.
For $T-r\ge10^{26}$,
\[
b_0(r)-1\ge10^{-20}-2\cdot10^{-24}-2e^{-100}.
\]
On $z\in[30,30+ch_0]$, $G(z)>1-10^{-24}$ and
$K(z)<3\cdot10^{-25}$. With $1\le\eta(a_*+s)\le2$, the integrand
in \mainequation{eq:strictgain} is greater than $9\cdot10^{-21}$.
The policy improvement is therefore strictly positive.
Choose time $t_*=1/2$. The stated lower bound on $T$
leaves at least $10^{26}$ units of remaining time throughout
the comparison, and $a_*\le t_*$.

The explicit inequalities can be verified with rational bounds.
For $n=60,100$, use
$e^{-n}<(\sum_{j=0}^{12}1/j!)^{-n}$ and
$e^{-1/40}>1-1/40$. These rational bounds give the inequalities
without relying on rounded decimals.
The resulting certified first-claim loss bounds are
$3.28\cdot10^{-7}$ above and $9.521484375\cdot10^{-6}$ below,
respectively; the policy-improvement lower bound exceeds
$4.9\cdot10^{-23}$. \end{proof}

\bibliographystyle{elsarticle-harv}
\bibliography{references}

@article{Avanzi2009,
 author={Avanzi, Benjamin},
 title={Strategies for dividend distribution: A review},
 journal={North American Actuarial Journal}, volume={13}, number={2},
 pages={217--251}, year={2009}, doi={10.1080/10920277.2009.10597549}
}

@book{Schmidli2008,
 author={Schmidli, Hanspeter}, title={Stochastic Control in Insurance},
 publisher={Springer}, address={London}, year={2008},
 doi={10.1007/978-1-84800-003-2}
}

@article{AlbrecherGarcia2023,
 author={Albrecher, Hansj{\"o}rg and Garcia Flores, Brandon},
 title={Optimal dividend bands revisited: A gradient-based method and evolutionary algorithms},
 journal={Scandinavian Actuarial Journal}, volume={2023}, number={8},
 pages={788--810}, year={2023}, doi={10.1080/03461238.2022.2147862}
}

@article{AvramPalmowskiPistorius2007,
 author={Avram, Florin and Palmowski, Zbigniew and Pistorius, Martijn R.},
 title={On the optimal dividend problem for a spectrally negative {L\'evy} process},
 journal={The Annals of Applied Probability}, volume={17}, number={1},
 pages={156--180}, year={2007}, doi={10.1214/105051606000000709}
}

@article{Loeffen2008,
 author={Loeffen, Ronnie L.},
 title={On optimality of the barrier strategy in de {Finetti}'s dividend problem for spectrally negative {L\'evy} processes},
 journal={The Annals of Applied Probability}, volume={18}, number={5},
 pages={1669--1680}, year={2008}, doi={10.1214/07-AAP504}
}

@article{AzcueMuler2005,
 author={Azcue, Pablo and Muler, Nora},
 title={Optimal reinsurance and dividend distribution policies in the {Cram\'er--Lundberg} model},
 journal={Mathematical Finance}, volume={15}, number={2},
 pages={261--308}, year={2005}, doi={10.1111/j.0960-1627.2005.00220.x}
}

@article{AlbrecherHartinger2006,
 author={Albrecher, Hansj{\"o}rg and Hartinger, J{\"u}rgen},
 title={On the non-optimality of horizontal barrier strategies in the {Sparre Andersen} model},
 journal={HERMIS}, volume={7}, pages={109--122}, year={2006}
}

@article{AlbrecherClaramuntMarmol2005,
 author={Albrecher, Hansj{\"o}rg and Claramunt, M. Merc\`e and M\'armol, Maite},
 title={On the distribution of dividend payments in a {Sparre Andersen} model with generalized {Erlang(n)} interclaim times},
 journal={Insurance: Mathematics and Economics}, volume={37}, number={2},
 pages={324--334}, year={2005}, doi={10.1016/j.insmatheco.2005.05.004}
}

@article{TianBaiGuo2020,
 author={Tian, Linlin and Bai, Lihua and Guo, Junyi},
 title={Optimal singular dividend problem under the {Sparre Andersen} model},
 journal={Journal of Optimization Theory and Applications}, volume={184}, number={2},
 pages={603--626}, year={2020}, doi={10.1007/s10957-019-01600-0}
}

@article{TianLiu2022,
 author={Tian, Linlin and Liu, Zhaoyang},
 title={Optimal dividend strategies in a renewal risk model with phase-type distributed interclaim times},
 journal={Applied Mathematics \& Optimization}, volume={85},
 note={Article 13}, year={2022}, doi={10.1007/s00245-022-09853-4}
}

@article{LiuLiuLiu2020,
 author={Liu, Yuying and Liu, Zhaoyang and Liu, Guoxin},
 title={Optimal dividend problems for {Sparre Andersen} risk model with bounded dividend rates},
 journal={Scandinavian Actuarial Journal}, volume={2020}, number={2},
 pages={128--151}, year={2020}, doi={10.1080/03461238.2019.1655475}
}

@article{DeAngelisEkstrom2017,
 author={De Angelis, Tiziano and Ekstr{\"o}m, Erik},
 title={The dividend problem with a finite horizon},
 journal={The Annals of Applied Probability}, volume={27}, number={6},
 pages={3525--3546}, year={2017}, doi={10.1214/17-AAP1286}
}

@book{AliprantisBorder2006,
 author={Aliprantis, Charalambos D. and Border, Kim C.},
 title={Infinite Dimensional Analysis: A Hitchhiker's Guide},
 edition={3}, publisher={Springer}, address={Berlin, Heidelberg},
 year={2006}, doi={10.1007/3-540-29587-9}
}

@article{MishuraSchmidli2012,
 author={Mishura, Yuliya and Schmidli, Hanspeter},
 title={Dividend barrier strategies in a renewal risk model with generalized {Erlang} interarrival times},
 journal={North American Actuarial Journal}, volume={16}, number={4},
 pages={493--512}, year={2012}, doi={10.1080/10920277.2012.10597645}
}

@article{LoeffenRenaud2010,
 author={Loeffen, Ronnie L. and Renaud, Jean-Fran{\c c}ois},
 title={De {Finetti}'s optimal dividends problem with an affine penalty function at ruin},
 journal={Insurance: Mathematics and Economics}, volume={46}, number={1},
 pages={98--108}, year={2010}, doi={10.1016/j.insmatheco.2009.09.006}
}

@article{KyprianouRiveroSong2010,
 author={Kyprianou, Andreas E. and Rivero, V{\'i}ctor and Song, Renming},
 title={Convexity and smoothness of scale functions and de {Finetti}'s control problem},
 journal={Journal of Theoretical Probability}, volume={23}, number={2},
 pages={547--564}, year={2010}, doi={10.1007/s10959-009-0220-z}
}

@inproceedings{DeFinetti1957,
 author={de Finetti, Bruno},
 title={Su un'impostazione alternativa della teoria collettiva del rischio},
 booktitle={Transactions of the XVth International Congress of Actuaries},
 volume={2}, pages={433--443}, address={New York}, year={1957}
}

@article{AsmussenTaksar1997,
 author={Asmussen, S{\o}ren and Taksar, Michael},
 title={Controlled diffusion models for optimal dividend pay-out},
 journal={Insurance: Mathematics and Economics}, volume={20}, number={1},
 pages={1--15}, year={1997}, doi={10.1016/S0167-6687(96)00017-0}
}
\end{document}